\documentclass[10pt, a4paper]{article}

\usepackage{amsmath,amsfonts,amssymb}
\usepackage{amsthm} % optional
\usepackage{hyperref}
\usepackage[nameinlink,capitalize]{cleveref}
\crefname{equation}{}{} 

\usepackage{enumerate}
\usepackage{float}
\usepackage[normalem]{ulem}
\usepackage{xcolor}

\usepackage{geometry}
\usepackage{bbm}
\usepackage{mathtools} % for vcentcolon

\usepackage{enumitem}
\setlist[enumerate]{label*=\alph*),ref=\alph*)}

\newcommand{\Asup}{\mathcal{A}_{\sup}}
\newcommand{\Anosup}{\mathcal{A}_{no\sup}}

\newcommand{\B}{\mathcal{B}}
\newcommand{\E}{\mathbb{E}}
\newcommand{\e}{\mathrm{e}}
\newcommand{\F}{\mathcal{F}}

\newcommand{\PP}{\mathbb{P}}
\newcommand{\R}{\mathbb{R}}
\newcommand{\N}{\mathbb{N}}
\newcommand{\one}{\mathbbm{1}}
\newcommand{\dd}{\mathrm{d}}
\newcommand{\define}{\vcentcolon =}
\newcommand{\EO}{\mathbb{E}_{\F_0}}
\newcommand{\PO}{\mathbb{P}_{\F_0}}
\newcommand{\pF}{_{p,\F_0}}

\newcommand{\xs}{X^{(n)}_{s^-}}
\newcommand{\g}{g(s, X^{(n)}_{\cdot \wedge k(n,s)},\xi)}
\newcommand{\h}{h(s, X^{(n)}_{\cdot \wedge k(n,s)})}
\newcommand{\nn}{^{(n)}}
\newcommand{\mm}{^{(m)}}
\newcommand{\tRnm}{\tau_R^{n,m}}
\newcommand{\cadlag}{\text{C\`adl\`ag}([-r, \infty);\R^d)}
\newcommand{\cadlagg}{c\`adl\`ag }
\newcommand{\kn}{\frac{k}{n}}
\newcommand{\levy}{L{\'e}vy }
\newcommand{\pred}{\mathcal{P}}
\newcommand{\Xt}{X_t}

\newcommand{\Xsm}{X_{s^-}}
\newcommand{\s}{_s}

\newcommand{\uu}{_u}
\newcommand{\sm}{_{s^-}}
\newcommand{\ti}{_t}
\newcommand{\dom}{\textrm{domain}(G^{-1})}

\newtheorem{theorem}{Theorem}[section]
\newtheorem{corollary}[theorem]{Corollary}
\newtheorem{lemma}[theorem]{Lemma}

\newtheorem{defi}[theorem]{Definition}

\theoremstyle{definition}
\newtheorem{example}[theorem]{Example}
\newtheorem{remark}[theorem]{Remark}

\newtheorem{hypo}[theorem]{Assumption}

\title{Concave and other generalizations of stochastic Gronwall inequalities}
\author{%
Sarah~Geiss\footnote{Technische Universit\"at Berlin, Germany. E-mail: \href{mailto:geiss@math.tu-berlin.de}{geiss@math.tu-berlin.de}}
\thanks{The author was supported by an 
Elsa-Neumann-Stipendium des Landes Berlin.} }

\date{April 20, 2023}   %% By default, LaTeX uses the current date

\begin{document}

\maketitle

\begin{abstract} 
We provide nonlinear generalizations of a class of stochastic Gronwall inequalities that have been studied by von Renesse and Scheutzow (2010), Scheutzow (2013), Xie and Zhang (2020) and Mehri and Scheutzow (2021). This class of stochastic Gronwall inequalities is a useful tool for SDEs.

More precisely, we study generalizations of the Bihari-LaSalle type. Whilst in a closely connected article by the author convex generalizations are studied, we investigate here concave and other generalizations. These types of estimates are useful to obtain existence and uniqueness of global solutions of path-dependent SDEs driven by \levy processes under one-sided non-Lipschitz monotonicity and coercivity assumptions. 
\end{abstract}

%\begin{thanks} TBA
%\end{thanks}

[0.25em]
\textbf{MSC2020 subject classifications:} 34K50, 60H10, 60G44, 60G51, 60J65

\bigskip

\tableofcontents

\section{Introduction}
In this article we provide generalizations of two types of stochastic Gronwall inequalities. Whereas in  the closely related article \cite{GeissConvex} convex generalizations are studied, we focus here on concave and more general generalizations. The estimates we obtain seem necessary to show existence, uniqueness and non-explosion of solutions to path-dependent SDEs under one-sided non-Lipschitz conditions. This section is structured as follows: First, we provide two simple examples which demonstrate where the generalizations we study can arise. Subsequently, we provide a short overview on the existing results in the literature and on the extensions we study in this paper.

\bigskip
\textbf{Path-dependent SDEs as motivation for stochastic Gronwall inequalities:} The stochastic Gronwall inequalities we generalize in this article were originally developed to study path-dependent SDEs \cite{RenesseScheutzow}: Consider the following path-dependent SDE driven by a Brownian motion $B$,
\begin{equation*}
dY_t = f(Y_{[-r,t]},t) \dd t + g(X_{[-r,t]},t) \dd B_t, \qquad Y_{[-r,0]} =y \in C([-r,0],\R^d),
\end{equation*}
where $X_{[-r,t]}$ denotes the path segment $\{X(u),u\in[-r,t]\}$ and $r>0$ is a constant. Uniqueness of global solutions is known under a one-sided Lipschitz condition
\begin{equation}\label{eq:intro-assumption}
2 \langle y(t)-z(t), f(y_{[-r,t]},t)-f(z_{[-r,t]},t)\rangle  + |g(y_{[-r,t]},t) - g(z_{[-r,t]},t)|_F^2 \leq K \sup_{s\in[-r,t]}|y(s)-z(s)|^2
\end{equation}
for all $y,z\in C([-r,\infty),\R^d)$, where $K>0$ is a constant.  Here $|\cdot|_F$ denotes the Frobenius norm. For stronger and more general results see von Renesse and Scheutzow \cite[Theorem 2.2]{RenesseScheutzow} and Mehri and Scheutzow \cite[Theorem 3.3]{MehriScheutzow}.
Note that \eqref{eq:intro-assumption} is weaker than assuming
\begin{equation*}
\begin{aligned}
2 \langle y(t)-z(t), f(y_{[-r,t]},t)-f(z_{[-r,t]},t)\rangle \leq   K \sup_{s\in[-r,t]}|y(s)-z(s)|^2  \\
\text{and} \quad |g(y_{[-r,t]},t) - g(z_{[-r,t]},t)|_F^2 \leq K \sup_{s\in[-r,t]}|y(s)-z(s)|^2
\end{aligned}
\end{equation*}
 as the scalar product term might be negative.  We sketch the proof ansatz for the uniqueness of global solutions: For two global solutions $Y$ and $Z$ with the same initial condition $y$ define the process $X_t \define |Y_t-Z_t|^2$. It suffices to show $\sup_{t\in[0,T]} X_t= 0$  a.s. for any $T>0$. By It\^o's formula and \eqref{eq:intro-assumption} we have
\begin{equation*}
\begin{aligned}
X_t & = \int_0^t 2 \langle Y_s-Z_s, f(Y_{[0,s]},s)-f(Z_{[0,s],s})\rangle  + |g(Y_{[0,s]},s) - g(Z_{[0,s]},s)|_F^2  \dd s  \\
&\qquad  + \int_0^t 2 \langle Y_s-Z_s,  g(Y_{[0,s]},s)-g(Z_{[0,s]},s) \dd B_s \rangle \\
& \leq \int_0^t K X_s^* \dd s + M_t \qquad \forall t\geq 0,
\end{aligned}
\end{equation*}
where  $X_s^*$ denotes the running supremum $\sup_{u\in[0,s]} X_u$ and $M$ is a continuous local martingale starting in $0$. Obtaining upper bounds for $X^*_t$ (i.e. showing it is $0$) using the BDG inequality and the Young inequality does not seem to be possible because then we would obtain the term $|g(y_{[0,t]},t) - g(z_{[0,t]},t)|_F^2$ due to the quadratic variation. Also taking expectations and then applying the deterministic Gronwall inequality to $t\mapsto \E[X_t]$ does not work, due to the running supremum $X^*_s$ in the integral. In \cite{RenesseScheutzow} bounds for $X$ are obtained instead by proving a stochastic Gronwall inequality: To this end it is exploited that $X$ is non-negative and $M$ is a continuous local martingale starting in $0$. In \cite{MehriScheutzow} these results are extended to \cadlagg local martingales.

Consider now a one-sided non-Lipschitzian monotonicity assumption: Replace in \eqref{eq:intro-assumption}  $K \sup_{s\in[0,t]}|y(s)-z(s)|^2$ by  $K \eta(\sup_{s\in[0,t]}|y(s)-z(s)|^2)$, where $\eta:[0,\infty) \to [0,\infty)$ is a non-decreasing function satisfying $\int_0^\varepsilon \frac{\dd u }{\eta(u)} = \infty$ for all $\varepsilon>0$ (Osgood condition).  To prove that uniqueness of global solutions also holds in this case, consider as before $X_t =|Y_t-Z_t|^2$. It would satisfy the following type of inequality
 \begin{equation*}
X_t \leq \int_0^t K \eta(X_s^*) \dd s + M_t \qquad \forall t\geq 0.
\end{equation*}
This is a special case of the inequalities (and applications) we study in this article.

\bigskip

\textbf{Application of stochastic Gronwall inequalities to other types of SDEs:}
Deterministic Gronwall inequalities are widely used to study SDEs. We provide a typical example when the deterministic Gronwall inequality is not applicable and instead stochastic generalizations of Gronwall inequalities appear useful to shorten calculations. Consider a non-path-dependent SDE driven by a Brownian motion $B$
\begin{equation*}
dY_t = b(Y_t,t) \dd t + \sigma(Y_t,t) \dd B_t, \quad Y_0 =y_0\in\R^d
\end{equation*}
with random coefficients $b$ and $\sigma$. Assume the following one-sided condition
\begin{equation}\label{eq:intro-asummption-0}
2 \langle y-z, b(y,t)-b(z,t)\rangle  + |\sigma(y,t) - \sigma(z,t)|_F^2 \leq K_t \eta(|y-z|^2) \qquad \forall y,z\in\R^d,\forall t\geq 0,
\end{equation}
where $\eta:[0,\infty) \to [0,\infty)$ is a non-decreasing function satisfying $\int_0^\varepsilon \frac{\dd u }{\eta(u)} = \infty$ for all $\varepsilon>0$, and $(K_t)_{t\geq 0}$ is an adapted non-negative stochastic process with suitable integrability assumptions.
It is known that under additional assumptions this  implies the uniqueness of global solutions. For more details in a more general setting see e.g. Lan and Wu \cite[Theorem 1.1]{WuLan} and the references therein. Let $Y$ and $Z$ be global solutions of this SDE and set  $X_t \define |Y_t-Z_t|^2$. Similarly as in the path-dependent case above, using It\^o's formula  to compute $X_t$ and applying assumption \eqref{eq:intro-asummption-0} gives
\begin{equation*}
\begin{aligned}
X_t \leq \int_0^t \eta(X_s) K_s \dd s + M_t \qquad \forall t\geq 0
\end{aligned}
\end{equation*}
for the local martingale $M_t \define \int_0^t 2 \langle Y_s - Z_s, \sigma(Y_s, s) - \sigma(Z_s, s)\rangle \dd s$. Hence, to prove uniqueness of solutions in this case, the typical argument to take expectations and then apply the deterministic Gronwall inequality to $t\mapsto \E[X_t]$ fails already in the linear case $\eta(x) = x$ if $(K_t)_{t\geq 0}$ is random. Using different techniques, uniqueness of solutions and further properties have been shown in Lan and Wu \cite[Theorem 1.1]{WuLan}.  Here, alternatively, generalizations of stochastic Gronwall inequalities could be applied and would shorten computations.

\bigskip

\textbf{Results in the literature on stochastic Gronwall inequalities:}  Now we summarize the existing literature on the two stochastic Gronwall inequalities we will generalize. The following  \emph{stochastic Gronwall inequality with supremum} is due to von Renesse and Scheutzow \cite[Lemma 5.4]{RenesseScheutzow} and was generalized by Mehri and Scheutzow \cite[Theorem 2.1]{MehriScheutzow}: Let $(\Xt)_{t\geq 0}$ be a non-negative stochastic process, that satisfies
\begin{equation}\label{eq:introsuplinear}
\Xt \leq  \int_{(0,t]} \Xsm^* \dd A\s + M\ti + H\ti \qquad \text{for all } t\geq 0,
\end{equation}
where $X^*\s= \sup_{u\in[0,s]}X\uu$  denotes the running supremum. Here, $M$ is a c\`adl\`ag local martingale that starts in $0$, and $H$ and $A$ are suitable non-decreasing stochastic processes. Then, for all $T>0$ and $p\in(0,1)$ there exists an explicit upper bound for $\E[\sup_{t\in[0,T]} X\ti^p]$ which does not depend on the local martingale $M$.

There is also a \emph{stochastic Gronwall inequality without supremum} which is closely connected to the previous inequality with supremum. This result is due to Scheutzow \cite[Theorem 4]{Scheutzow} and was generalized by Xie and Zhang \cite[Lemma 3.7]{XieZhang} to c\`adl\`ag local martingales: If we assume instead of \eqref{eq:introsuplinear} the slightly stronger assumption that 
\begin{equation}\label{eq:intronosuplinear}
X\ti \leq  \int_{(0,t]} X\sm \dd A\s + M\ti + H\ti \qquad \text{for all } t\geq 0,
\end{equation}
sharper bounds can be obtained.

Both previously mentioned inequalities are useful tools for SDEs. The stochastic Gronwall inequality with supremum 
is applied to study SDEs with memory, see for  example  \cite{AppSup1}, \cite{AppSup4},  \cite{AppSup3}, \cite{AppSup2}, \cite{HutzenthalerNguyen}, \cite{MehriScheutzow}, \cite{AppSup5} and \cite{RenesseScheutzow}. 
The stochastic Gronwall inequality without supremum is applied to study various SDEs without memory, see e.g.
 \cite{AppNoSup4}, \cite{AppNoSup3}, \cite{AppNoSup8}, \cite{AppNoSup10}, \cite{AppNoSup11},
 \cite{AppNoSup1}, \cite{AppNoSup7}, \cite{AppNoSup5}, \cite{AppNoSup2}, \cite{AppNoSup6} and
 \cite{AppNoSup9}. 
 
Further stochastic Gronwall inequalities with a somewhat different setting have been studied by Glatt-Holtz and Ziane \cite[Lemma 5.3]{GlattZiane} and Agresti and Veraar \cite[Lemma A.1]{AgrestiVeraar}.
\bigskip

\bigskip \bigskip
\textbf{Results of this article on generalizations of stochastic Gronwall inequalities:} In this paper, we study the following nonlinear generalization of the above mentioned stochastic Gronwall inequalities: We replace the assumptions \eqref{eq:introsuplinear} and \eqref{eq:intronosuplinear} by
\begin{equation}\label{eq:introsup}
X\ti \leq  \int_{(0,t]} \eta(X^*\sm) d A\s + M\ti + H\ti \qquad \text{for all } t\geq 0,
\end{equation}       
and
\begin{equation}\label{eq:intronosup}
X\ti \leq  \int_{(0,t]} \eta(X\sm) d A\s + M\ti + H\ti \qquad \text{for all } t\geq 0,
\end{equation}
respectively, where $\eta:[0,\infty) \to [0,\infty)$ is a non-decreasing function. For continuous martingales and $\eta$ that satisfy $\int_1^\infty \frac{\dd u }{\eta(u)} = +\infty$ inequality \eqref{eq:introsup} has been studied by von Renesse and Scheutzow \cite[Lemma 5.1]{RenesseScheutzow} in the context of global solutions of stochastic functional differential equations, but no explicit upper bounds for $\E[\sup_{t\in[0,T]} X\ti^p]$ were derived. For \emph{convex} $\eta$ these inequalities are studied in \cite{GeissConvex}. Also other types of nonlinear extensions of the stochastic Gronwall inequalities have been studied, see \cref{subsec:results} for more details.

For the deterministic case (i.e. $M\equiv 0$) the Bihari-LaSalle inequality provides an upper bound for $X\ti$, see \cite{Bihari}, \cite{LaSalle}. In \cite{Bihari} and \cite{LaSalle} the integrator is of the form $A(t) = \int_0^t \varphi(s) ds$, but we will need a version where the integrator $A$ is only c\`adl\`ag. As this version is difficult to find in the literature, we provide a short proof in the appendix, which goes along the lines of the standard proof.

\begin{lemma}[Deterministic Bihari-LaSalle inequality]\label{lemma:DetBihari}
Let $H>0$ be a constant and let $x:[0,\infty)\to[0,\infty)$ be a c\`adl\`ag function and  $A:[0,\infty)\to[0,\infty)$ a non-decreasing c\`adl\`ag function with $A(0)=0$. Let $\eta:[0,\infty) \to [0,\infty)$ be a left-continuous non-decreasing function satisfying $\eta(u)>0$ for $u>0$. Define $G(v) \define \int_c^v \frac{\dd u}{\eta(u)}$ for some constant $c>0$.
Assume that the function $x$ satisfies
\begin{equation}\label{eq:detBihari}
x(t) \leq \int_{(0,t]} \eta(x(s^-)) \dd A(s) + H \qquad \forall t\in[0,T]
\end{equation}
for some $T>0$. If $G(H)+A(T) \in \textrm{domain}(G^{-1})$, then the following inequality holds true:
\begin{equation*}
x(t) \leq G^{-1}(G(H) + A(t)) \qquad \forall t\in[0,T].
\end{equation*}
\end{lemma}

\smallskip 

\noindent Note that the upper bound on $x$ does not depend on the choice of the constant $c>0$ used in the definition of $G$. We obtain the well-known Gronwall inequality by choosing $c=1$ and $\eta(u) \equiv u$ in \cref{lemma:DetBihari} so that $G(u) \equiv \log(u)$, which implies the upper bound
\begin{equation*}
x(t)\leq H \exp(A(t)).
\end{equation*}

In the present article we show for concave $\eta$ a stochastic version of the Bihari-LaSalle inequality which provides upper bounds similar to that of the deterministic Bihari-LaSalle inequality. Let $(\Omega, \F, \PP, (\F_t)_{t\geq 0}$ an underlying filtered probability space satisfying the usual conditions. The following theorem is a part of \cref{thm:stochBihariConcave} (noting that $x\mapsto G^{-1}(G(x^{1/p})+\beta a)^p$ is concave under the assumptions of the theorem, i.e. $\E[\cdot \mid \F_0]$ may be replaced by $\E[\cdot]$ by Jensen's inequality).
\begin{theorem}[Concave stochastic Bihari-LaSalle inequality with supremum]
Let $p\in(0,1)$ and assume that 
\begin{itemize}
\item  $(X_t)_{t\geq 0}$ is an adapted non-negative right-continuous process,
\item  $(A_t)_{t\geq 0}$ is a deterministic non-decreasing c\`adl\`ag function with $A_0=0$, 
\item  $(H_t)_{t\geq 0}$ is an predictable non-negative non-decreasing c\`adl\`ag process,
\item  $(M_t)_{t\geq 0}$ is a c\`adl\`ag local martingale with $M_0=0$,
\item  $\eta:[0,\infty) \to [0,\infty)$ is a continuous non-decreasing function with $\eta(x)>0$ for all $x>0$ such that $\eta_p:(0,\infty) \mapsto [0,\infty), \,\, x\mapsto \tfrac{p}{1-p} \eta(x^{1/p})x^{1-1/p}$ is concave and non-decreasing.
\end{itemize}
Moreover, assume that 
\begin{equation*}
X\ti \leq \int_{(0,t]} \eta(X^*\sm) \dd A\s + M\ti + H\ti \qquad \PP\text{-a.s}.
\end{equation*}
holds true for all $t\geq 0$. Then, for all $T>0$ such that $\E[H_T^p]<\infty$ the following inequality holds
\begin{equation*}
\E\left[\sup_{t\in[0,T]}X^p_t\right]^{1/p} \leq G^{-1}\big(G(\alpha_1 \alpha_2 \E[H_T^p]^{1/p})+ \beta A_T\big),
\end{equation*}
where we use the notation $G(v) \define \int_1^v \frac{\dd u}{\eta(u)}$. The constants $\beta = (1-p)^{-1}$ and $\alpha_1 \alpha_2 = (1-p)^{-1/p}p^{-1}$  are sharp.
\end{theorem}
The assumption that $\eta_p$ is concave is slightly stronger than the assumption that $\eta$ is concave. By a time change argument this theorem also implies estimates for random integrators $A$, see \cref{cor:randomA}.

Under the stronger assumption \eqref{eq:intronosup} we obtain for concave $\eta$ and random $A$ also the following weak $L^1$ estimate: For all $t\geq 0$, $R>0$, $w> 0$, $u> 0$ we have
\begin{equation}\label{eq:intro-veraar}
\PP[X^*_t > u] \leq 
\frac{G^{-1}(G(\E[H_t\wedge w]) + R)}{u} + \PP[H_t \geq w] + \PP[A_t >  R],
\end{equation}
see \cref{cor:randomAnoSup}. The benefit of this estimate is that it does not require integrability assumptions on $H$ and $A$. It can be shown that under the weaker assumption \eqref{eq:introsup} the weak $L^1$ norm is in general not finite, in particular 
\eqref{eq:intro-veraar} fails to hold true in this case. Note that in the special case $\eta(x)=x$ the estimate  \eqref{eq:intro-veraar} simplifies to
\begin{equation}\label{eq:intro-veraar-linear}
\PP[X^*_t > u] \leq 
\frac{\e^R }{u}\E[H_t\wedge w] + \PP[H_t \geq w] + \PP[A_t >  R],
\end{equation}
for more details on the estimate \eqref{eq:intro-veraar-linear} and the linear case $\eta(x)= x$ see \cite[Section 5]{GeissConvex}. The inequality \eqref{eq:intro-veraar-linear} is applied by Agresti and Veraar  \cite[Lemma 5.3, Lemma 6.7]{AgrestiVeraar2} to prove global well-posedness for reaction-diffusion systems with transport noise.

For the application to path-dependent SDEs the cases $\int_{0}^{1} \frac{\dd u}{\eta(u)} = +\infty$ and $\int_1^\infty \frac{\dd u}{\eta(u)} = +\infty$ (without necessarily imposing concavity or convexity assumptions on $\eta$) are especially useful.  Whilst we do not obtain explicit upper bounds for the quantity $\E[\sup_{t\in[0,T]}X^p_t]$ without  concavity or convexity assumptions, we do observe that 
 $(X_t)_{t\in[0,T]}$ still behaves similarly w.r.t the process $(H_t)_{t\geq 0}$ as the function $x(t), t\in[0,T]$ w.r.t. the constant $H$ in the deterministic Bihari-LaSalle inequality, for details see \cref{thm:BihariRandomA}.  These properties are sufficient for our applications: We apply these results to prove existence and uniqueness of global solutions to path-dependent SDEs driven by L\'evy processes under one-sided non-Lipschitzian monotonicity and coercivity assumptions.

%%%%%%%%%%%%%%%%%%%%%%%%%%%%%%%%%%%%%%%%%%%%%%%%%%%%%%%%%%%%%%%%%%%
%%%%%%%%%%%%%%%%%%%%%%%%%%%%%%%%%%%%%%%%%%%%%%%%%%%%%%%%%%%%%%%%%%%
%%%%%%%%%%%%%%%%%%%%%%%%%%%%%%%%%%%%%%%%%%%%%%%%%%%%%%%%%%%%%%%%%%%

\section{Notation, assumptions and overview}
We assume that all processes are defined on an underlying filtered probability space $(\Omega, \F, \PP, (\F_t)_{t\geq 0})$ satisfying the usual conditions.

\subsection{Notation and constants} \label{subsec:constants}
\textbf{Constants:} For $p\in(0,1)$ define the following constants:
\begin{equation}\label{eq:constants}
\beta = (1-p)^{-1},\quad   \alpha_1 = 
    (1-p)^{-1/p}, \quad  \alpha_2 = p^{-1}.
\end{equation}
\textbf{Quasinorms:} We denote by $|\cdot|$  the Euclidean norm and $|\cdot|_F$ the Frobenius norm. Let $Y$ be a random variable. We use for $p\in(0,1]$ the notation (if well-defined)
\begin{equation*}
\EO[Y] \define \E[\,Y\mid \F_0], \qquad \qquad 
\|Y\|_p \define \E[\,|Y|^p]^{1/p}, \qquad  \qquad \|Y\|_{p,\F_0} \define \E[\,|Y|^p\mid \F_0]^{1/p}.
\end{equation*}

\noindent
\textbf{Running supremum:} Let $X$ be a non-negative stochastic process with right-continuous
%c\`adl\`ag 
paths. We use the following notation for the running supremum and its left limits:
\begin{equation*}
\begin{aligned}
X^*\ti  & \define \sup_{0\leq s \leq t} X\s \qquad \forall t\geq 0, \\
X^*_{t^-} & \define \lim_{s\nearrow t} X^*\s = \sup_{s< t}X\s \qquad \forall t>0.
\end{aligned}
\end{equation*}
As usual, we set $X^*_{0^-} \define X_0$. If $X$ is c\`adl\`ag, then also $X^*_{t^-} = \sup_{s\leq t} X\sm$ holds true. If $X$ is only right-continuous then $X^*_t$ and $X^*_{t^-}$ take values in $[0,+\infty]$.\\

\noindent \textbf{Functions:} Let $\eta:[0,\infty) \to [0,\infty)$ be a non-decreasing function with $\eta(u) > 0$ for $u>0$. Define the following functions for $p\in(0,1)$.
%\begin{equation}
\begin{align}
G(x) & \define \int_c^x \frac{\dd u}{\eta(u)} \qquad & \forall x\in[0,\infty), \label{eq:defG}\\
\eta_p(x) & \define \frac{p}{1-p}\eta(x^{1/p})x^{1-1/p} \qquad   &\forall x\in(0,\infty), \label{eq:defEta_p} \\
\tilde G_p(x) & \define \int_{c^{p}}^{x} \frac{\dd u}{\eta_p(u)} \qquad  & \forall x\in(0,\infty).\label{eq:defTildeG}
\end{align}
%\end{equation*}
The functions have the following properties:
\begin{itemize}
\item The function $G$ satisfies $G(0)\in[-\infty,0)$. Moreover, $G$ is increasing and concave. In particular, it has a well-defined increasing inverse $G^{-1}:\textrm{range}(G)\cap(-\infty,\infty) \mapsto [0,\infty)$. If $G(0) = -\infty$, then we set $G^{-1}(G(0)+a)\define 0$ for $a\in(-\infty,\infty)$. If $\eta$ is continuous, then $G$ is continuously differentiable on $(0,\infty)$. 
\item The functions $G$ and $\tilde{G}_p$ satisfy for all $x\in (0,\infty)$
\begin{equation}\label{eq:G_umrechnen}
\begin{aligned}
\tilde G_p(x) & = \frac{1-p}{p}\int_{c^p}^{x}\frac{\dd u}{ \eta(u^{1/p}) u^{1-1/p}} =  (1-p)\int_{c}^{x^{1/p}}\frac{\dd v}{ \eta(v)}  &= (1-p)G(x^{1/p})
\end{aligned}
\end{equation}
and for all $x\in\textrm{domain}(\tilde{G}^{-1}_p)$
\begin{equation}\label{eq:G_umrechnen2}
 \tilde{G}^{-1}_p(x) = \big(G^{-1}\big(\tfrac{x}{1-p} \big)\big)^p.
\end{equation}
\end{itemize}

\subsection{Assumptions}\label{subsec:assumptions}
We study the following two cases:
\begin{defi}[Assumption $\mathcal{A}_{\sup}$] \label{def:sup}
Let
\begin{itemize}
\item  $(X_t)_{t\geq 0}$ be an adapted non-negative right-continuous process,
\item  $(A_t)_{t\geq 0}$ be a predictable non-decreasing c\`adl\`ag process with $A_0=0$, 
\item  $(H_t)_{t\geq 0}$ be an adapted non-negative non-decreasing c\`adl\`ag process,
\item  $(M_t)_{t\geq 0}$ be a c\`adl\`ag local martingale with $M_0=0$,
\item  $\eta:[0,\infty) \to [0,\infty)$ be a continuous non-decreasing function with $\eta(x)>0$ for all $x>0$.
\end{itemize}
We say the processes $X$, $A$, $H$, $M$ satisfy $\mathcal{A}_{\sup}$ if they satisfy the inequality below for all $t\geq 0$:
\begin{equation}\label{eq:GronwallAssumption}
X\ti \leq \int_{(0,t]} \eta(X^*\sm) \dd A\s + M\ti + H\ti \qquad \PP\text{-a.s}.
\end{equation}
\end{defi}

\noindent The following assumption is slightly stronger:
\begin{defi}[Assumption $\mathcal{A}_{no\sup}$]\label{def:nosup}
 Under the same assumptions on the processes as in the previous definition, we say that the processes satisfy $\mathcal{A}_{no\sup}$ if in addition $X$ has left limits and the processes satisfy the following inequality for all $t\geq 0$:

\begin{equation}\label{eq:GronwallAssumptionNoSup}
X\ti \leq \int_{(0,t]} \eta(X\sm) \dd A\s + M\ti + H\ti \qquad \PP\text{-a.s}.
\end{equation}
\end{defi}

We also use the two definitions above for processes defined on a finite time interval $[0,T]$ and correspondingly adapt the definition in this case.

\subsection{Overview of the results}\label{subsec:results}
\cref{table:Gronwall} summarizes the results obtained by von Renesse and Scheutzow
\cite[Lemma 5.4]{RenesseScheutzow}, Scheutzow \cite[Theorem 4]{Scheutzow}, Xie and Zhang \cite[Lemma 3.7]{XieZhang}, Mehri and Scheutzow \cite[Theorem 2.1]{MehriScheutzow} and the author \cite[Corollary 5.2, Corollary 5.4]{GeissConvex} on stochastic Gronwall inequalities.
For a detailed overview on these results, see e.g. \cite[Section 5]{GeissConvex}.
\cref{table:BihariLaSalle} summarizes the stochastic Bihari-LaSalle inequalities obtained in this article, and for the convenience of the reader, also the results obtained in \cite{GeissConvex}. 

The constant $\alpha_1\alpha_2$ is the sharp constant from Lenglart's inequality, $\alpha_2$ is the sharp constant from a monotone version of Lenglart's inequality \cite{GeissScheutzow}. The constants $\alpha_1$, $\alpha_1\alpha_2$ and $\beta$  appearing in the stochastic Gronwall inequalities are sharp, see \cite[Theorem 3.5, Theorem 3.7]{GeissConvex}.

\begin{table}[H]
\bgroup
\scriptsize
\def\arraystretch{2}
    \begin{tabular}{|p{2.2cm}|p{5.9cm}|p{6.1cm}|}
         \hline
             &  Assumption $\Anosup, \quad \eta(x) \equiv x$ \newline (Special case of $\Asup$)  & Assumption $\Asup, \quad \eta(x) \equiv x$ \\ \hline            
             $A$ deterministic,   \newline $p\in(0,1)$  
             & % without sup
              $\blacktriangleright$ $H$ predictable or  
              $\Delta M \geq 0$: \newline
              $\|X^*_T\|_p \leq \alpha_1 \alpha_2 \|H_T\|_p 
              \e^{A_T} $ 
              \newline and \newline
          $\PP[X_T^* > u] \leq \frac{e^{A_T}}{u} \E[H_T\wedge w] 
          + \PP[H_T \geq w]$
               \newline\newline
              $\blacktriangleright$ $\E[H_T] <\infty$: \newline
              $\|X^*_T\|_p \leq \alpha_1 \|H_T\|_1 \e^{A_T}$ 
              \newline and \newline
          $\PP[X_T^* > u] \leq \frac{e^{A_T}}{u} \E[H_T]$
             & %with sup
              $\blacktriangleright$ $H$ predictable or  
              $\Delta M \geq 0$: \newline
              $\|X^*_T\|_p \leq \alpha_1 \alpha_2 \|H_T\|_p \e^{\beta A_T} $ 
               \newline\newline
              $\blacktriangleright$ $\E[H_T] <\infty$: \newline
              $\|X^*_T\|_p \leq \alpha_1 \|H_T\|_1 \e^{\beta A_T}$ \\
            \hline
            $A$ random, \newline $0<q<p<1$
              & %without sup
              $\blacktriangleright$ $H$ predictable or 
               $\Delta M \geq 0$ \newline
              $\|X^*_T\|_q \leq \alpha_1 \alpha_2 
              \|H_T\|_p \|e^{A_T}\|_{qp/(p-q)} $ \newline
               and 
               \newline $\PP[X_T^* > u] $ \newline $\leq  \frac{e^{R}}{u}\E[H_T \wedge w ] + 
               \PP[H_T\geq w] + \PP[A_T> R]$ 
               \newline\newline\newline
              $\blacktriangleright$ 
              $\E[H_T] <\infty$: \newline
              $\|X^*_T\|_q \leq \alpha_1 \|H_T\|_1 \|e^{A_T}\|_{qp/(p-q)}$ 
              \newline
              and 
              \newline $\PP[X_T^* > u] \leq  \frac{e^{R}}{u}\E[H_T] + \PP[A_T> R]$
              & %with sup
              $\blacktriangleright$ $H$ predictable or  $\Delta M \geq 0$: \newline
              $\|X^*_T\|_q \leq \alpha_1\alpha_2 \|H_T\|_p 
              \|e^{\beta A_T}\|_{qp/(p-q)} $
              \newline\newline
              $\blacktriangleright$ $\E[H_T] <\infty$:  \newline
              $\|X^*_T\|_q \leq \alpha_1 \|H_T\|_1 \|e^{\beta A_T}\|_{qp/(p-q)}$  \\
     \hline
         Constants & \multicolumn{2}{c|}{$\beta = (1-p)^{-1},\quad   \alpha_1 = 
    (1-p)^{-1/p}, \quad  \alpha_2 = p^{-1}$}  \\
      \hline
      Notation  & \multicolumn{2}{c|}{$\|Y\|_p \define \E[|Y|^p]^{1/p}$ for random variables $Y$}  \\
       \hline 
    \end{tabular}
    \egroup
\caption{Summary of stochastic Gronwall inequalities (i.e. $\eta(x)\equiv x$) known in the literature}
\label{table:Gronwall}
\end{table}

\begin{table}[H]
\bgroup
\scriptsize
\def\arraystretch{2}
    \begin{tabular}{|p{2.07cm}|p{6.05cm}|p{6.23cm}|}
         \hline
             &  Assumption $\Anosup$\newline (Special case of $\Asup$) & Assumption $\Asup$  
             \\ 
             \hline
             $A$ random & 
             & \raggedright 
             $\blacktriangleright$  For $\int_1^\infty
             \frac{\dd u}{\eta(u)} = +\infty$ and continuous
             $M$ \newline see von Renesse and Scheutzow
             \cite[Lemma 5.1]{RenesseScheutzow}
             \newline \newline
             $\blacktriangleright$ If $\lim_{x\to\infty}\eta(x) = +\infty$, $\E[H_T]<\infty$ and 
             $\E[A^p_T]<\infty$:
             \newline
             $\|G(c+X^*_T)\|_p \leq 
             \alpha_1 \alpha_2 \|A_T +  G(c+\E[H_T])\|_p$ 
             \newline for $p\in(0,1)$, $G(x) \define \int_c^x 
             \frac{\dd u}{\eta(u)}$, $c>0$             
             \newline See \cref{thm:BihariRandomA}
             \newline\newline
             $\blacktriangleright$ If $\int_{0+}
             \frac{\dd u}{\eta(u)} = +\infty$ and 
             $\E[H_T]=0$: \newline
             $X_T^*=0 \quad \PP$-a.s. \newline
             See \cref{thm:BihariRandomA}
             \tabularnewline
             \hline            
             $\eta$ concave, \newline  $A$ deterministic, \newline $p\in(0,1)$    
             & % without sup
              $\blacktriangleright$ $H$ predictable or  
              $\Delta M \geq 0$: \newline
              $\|X^*_T\|_p \leq \alpha_1 G^{-1}(G(\alpha_2 
              \|H_T\|_p)+A_T) $ 
              \newline See \cref{thm:stochBihariConcave}
               \newline\newline
              $\blacktriangleright$ $\E[H_T] <\infty$: \newline
              $\|X^*_T\|_p \leq \alpha_1 G^{-1}(G(\|H_T\|_1)+  A_T)$
              \newline and \newline
              $\PP[X_T^*>u] \leq \frac{G^{-1}(G(\E[H_T]) 
              + A_T)}{u}$
              \newline See \cref{thm:stochBihariConcave} 
             & %with sup
             $\blacktriangleright$ $H$ predictable or  $\Delta M \geq 0$: \newline
              $\|X^*_T\|_p \leq G^{-1}(G(\alpha_1\alpha_2\|H_T\|_p)+\beta A_T) $ 
              \newline See \cref{thm:stochBihariConcave} \newline\newline
              $\blacktriangleright$ $\E[H_T] <\infty$: \newline
              $\|X^*_T\|_p \leq G^{-1}(G(\alpha_1 \|H_T\|_1)+\beta A_T)$ \newline See \cref{thm:stochBihariConcave}
             \\ 
             \hline
             $\eta$ concave, \newline  $A$ random, 
             \newline $p\in(0,1)$ 
             \newline $p<q$
             & % A random, without sup
             
              $\blacktriangleright$ $\E[H_T] <\infty$: \newline
              $\PP[X^*_T > u] \leq \frac{G^{-1}(G(\E[H_T]) + R)}{u} 
              + \PP[A_T >  R]$              
             \newline See \cref{cor:randomAnoSup} 
             & % A random, with sup 
             \raggedright 
              $\blacktriangleright$ $\E[H_T] <\infty$: \newline
             $\|X^*_t\|_{p/2} $ \newline 
             $\leq  C\|H_t\|_{1}^{(p-q)/p}\E[\big(G^{-1}
             (G(\alpha_1 \|H_t\|_{1}) 
             + \beta A_t)\big)^{q}]^{1/p}$
             \newline  
            for $C = (\frac{q}{q-p})^{1/p} \alpha_1^{(p-q)/p}$
              \newline See \cref{cor:randomA}
             \tabularnewline
             \hline
             $\eta$ convex,  \newline $A$ predictable, 
             \newline $p\in(0,1)$
             &
             % without sup
             $\blacktriangleright$ $H$ predictable or  
              $\Delta M \geq 0$: \newline
               $\|\sup_{t\in[0,T]} G^{-1}(G(X_t) - A_t)\|_p$ 
               \newline $\phantom{aaa }\leq
              \alpha_1\alpha_2 \|H_T\|_p$ \newline
              See  \cite[Theorem 3.1]{GeissConvex}
              \newline \newline
             $\blacktriangleright$  $\E[H_T] <\infty$: \newline 
            $\|\sup_{t\in[0,T]} G^{-1}(G(X_t) - A_t)\|_p  \leq
              \alpha_1 \|H_T\|_1$
              \newline and \newline
              $\PP[X_T^*>u] \leq \frac{\E[H_T]}{G^{-1}(G(u)-R)} 
              + \PP[A_T>R]$
              \newline
              See  \cite[Theorem 3.1]{GeissConvex}
             & $\blacktriangleright$ $H$ predictable or  
              $\Delta M \geq 0$: \newline
               $\|G^{-1}(G(X_T^*) - \beta A_T)\|_p  \leq
              \alpha_1\alpha_2 \|H_T\|_p$ \newline
              See  \cite[Theorem 3.1]{GeissConvex}
              \newline \newline
             $\blacktriangleright$  $\E[H_T] <\infty$: \newline 
            $\|G^{-1}(G(X_T^*) - \beta A_T)\|_p  \leq
              \alpha_1 \|H_T\|_1$
              \newline
              See  \cite[Theorem 3.1]{GeissConvex}\\
         \hline
         Constants & \multicolumn{2}{c|}{$\beta = (1-p)^{-1},\quad   \alpha_1 = 
    (1-p)^{-1/p}, \quad  \alpha_2 = p^{-1}$}  \\
     \hline
     Notation  & \multicolumn{2}{c|}{$\|Y\|_p \define \E[|Y|^p]^{1/p}$ for random variables $Y$, $\quad G(x) \define \int_c^x \frac{\dd u}{\eta(u)}$ $\quad \forall x\geq c_0$}  \\
       \hline 
    \end{tabular}
    \egroup
\caption{Summary of stochastic Bihari-LaSalle inequalities}
\label{table:BihariLaSalle}
\end{table}

Also other nonlinear extensions of the stochastic Gronwall inequalities have been studied. Makasu \cite[Theorem 2.2]{Makasu1} and Le and Ling \cite[Lemma 3.8]{LeLing} studied a non-linear generalization where in the assumption the term $\int_{(0,t]} X\sm \dd A\s$  is replaced by $\big(\int_{(0,t]} (X_s)^\theta \dd A\s\big)^{1/\theta}$ for $\theta >0$ and obtained estimates for $\E[\sup_{t\in[0,T]}X^p_t]$.
 
 A further nonlinear extension of the stochastic Gronwall inequalities has been studied by Mekki, Nieto and Ouahab \cite[Theorem 2.4]{Nieto}: For continuous local martingales a stochastic Henry Gronwall's inequality 
with upper bounds that do not depend on the local martingale $M$ can be proven.

In addition, also extensions of the stochastic Gronwall inequalities have been studied, where the upper bounds depend on the quadratic variation of the martingale $M$, see e.g. Makasu \cite{Makasu1}, Makasu \cite{Makasu2} and Mekki, Nieto and Ouahab \cite{Nieto}. In the present paper, we focus on bounds which do not depend on the local martingale $M$.  Furthermore,  Hudde, Hutzenthaler and Mazzonetto \cite{Hudde} have extended the stochastic Gronwall inequality without supremum to the setting of It\^o processes which satisfy a suitable one-sided affine-linear growth condition.

\section{Main results: Stochastic Bihari-LaSalle inequalities}

In this section we provide nonlinear generalizations of the stochastic Gronwall inequalities mentioned in the introduction and \cref{table:Gronwall}. For the definition of the constants see \eqref{eq:constants}. 

\subsection{Estimates for concave $\eta$ and deterministic $A$}

\begin{theorem}[A sharp stochastic Bihari-LaSalle inequality for concave $\eta$] \label{thm:stochBihariConcave}
Let $(X,\,A,\,H,\,M)$ and $\eta$ satisfy \nameref{def:sup} or \nameref{def:nosup}. Moreover, assume that $A$ is deterministic and $p\in(0,1)$. For $\eta$ we use the notation
\begin{equation*}
\eta_p:(0,\infty) \mapsto [0,\infty), \qquad x\mapsto \tfrac{p}{1-p} \eta(x^{1/p})x^{1-1/p}.
\end{equation*}

\begin{enumerate}
\item\label{item:concaveSup} Let \nameref{def:sup} hold and assume that $\eta_p$ is concave and non-decreasing. Then, the following assertions hold for all $t\geq 0$:
\begin{equation*}
\|X^*\ti\|_{p,\F_0}  \leq 
\begin{cases}
G^{-1} \bigg(G\big( \alpha_1 \alpha_2 \|H\ti \|_{p,\F_0}\big) 
+ \beta A\ti \bigg) & \text{ if } \E[H^p\ti] < \infty \text{ and } H \text{ is predictable}, \\[1em]
G^{-1} \bigg(G\big( \alpha_1 \alpha_2 \|H\ti \|_{p,\F_0}\big) 
+ \beta A\ti \bigg) & \text{ if } \E[H^p\ti] < \infty \text{ and }\Delta M \geq 0, \\[1em]
G^{-1} \bigg(G\big( \alpha_1 \|H\ti\|_{1,\F_0} \big) 
+ \beta A\ti\bigg) & \text{ if } \E[H\ti] < \infty,
\end{cases}
\end{equation*}
where $\beta\define(1-p)^{-1}$, $\alpha_1 \define(1-p)^{-1/p}$ and $\alpha_2 \define p^{-1}$. We use the notation $(\Delta M_t)_{t\geq 0}  \define (M_t - M_{t^-})_{t\geq 0}$ and $\|Y\|_{p,\F_0} \define \E[\,|Y|^p\mid \F_0]^{1/p}$ for random variables $Y$.

The constants $\alpha_1$, $\alpha_1\alpha_2$ and $\beta$ are sharp, they are already sharp when $\eta(x) \equiv x$.

\item\label{item:concaveNoSup} If \nameref{def:nosup} holds, the following assertions are true for all $t\geq 0$:

\begin{equation}\label{eq:claimConcaveNoSup-p}
\|X^*\ti\|_{p,\F_0}  \leq 
\begin{cases}
\alpha_1  G^{-1} \bigg(G\big(\alpha_2 \|H\ti \|_{p,\F_0}\big) 
+ A\ti \bigg)  &\text{ if } \E[H^p\ti] < \infty,\,\, H \text{ is predictable}, \\
& \text{and } \eta_p \text{ is concave and non-decreasing} \\[1em]
\alpha_1  G^{-1} \bigg(G\big(\alpha_2 \|H\ti \|_{p,\F_0}\big) 
+ A\ti \bigg) &\text{ if } \E[H^p\ti] < \infty,\,\,\Delta M \geq 0, \\
& \text{and } \eta_p \text{ is concave and non-decreasing}, \\[1em]
\leq \alpha_1 G^{-1} \bigg(G\big(  \|H\ti\|_{1,\F_0} \big) 
+ A\ti \bigg) &\text{ if } \E[H\ti] < \infty \text{ and } \eta \text{ is concave}.
\end{cases}
\end{equation}
The constant $\alpha_1$ is sharp. If $\eta(x) \equiv x$, then $\alpha_1 \alpha_2$ appears in the estimates for $\E[H^p\ti]<\infty$ and the constant is sharp.

Under \nameref{def:nosup} also the following weak $L^1$ estimates are true for all $t\geq 0$, $u> 0$, $w>0$:
\begin{equation}\label{eq:claimConcaveNoSup-weak}
\PO[X^*_t > u] \leq 
\begin{cases}
\frac{G^{-1}(G(\EO[H_t\wedge w]) + A_t)}{u} + \PO[H_t \geq w]
 &\text{ if } H \text{ is predictable and } \eta \text{ is concave,} \\[1em]
\frac{G^{-1}(G(\EO[H_t]) + A_t)}{u}
 &\text{ if } \E[H\ti] < \infty \text{ and } \eta \text{ is concave}.
\end{cases}
\end{equation}

\end{enumerate}
\end{theorem}

\begin{remark}
Assume that $(X,\,A,\,H,\,M)$ and $\eta$ satisfy \nameref{def:nosup} and let $\eta$ be concave. Then, also for the case $\Delta M \geq 0$ it can be shown that  $t\geq 0$, $u>0$, $w>0$ we have
\begin{equation*}
\PO[X^*_t > u] \leq 
\frac{G^{-1}(G(\EO[H_t\wedge w]) + A_t)}{u} + \PO[H_t \geq w].
\end{equation*}
We prefer not to prove it to keep the manuscript at a reasonable length.
\end{remark}

\begin{remark}\label{rmk:noWeak}
If we assume that $(X,\,A,\,H,\,M)$ and $\eta$ satisfy \nameref{def:sup}  (which is a weaker assumption than  \nameref{def:nosup}), \eqref{eq:claimConcaveNoSup-weak} does not hold in general. See \cite[Theorem 3.8]{GeissConvex} for a counterexample and further details.
\end{remark}

\begin{remark}\label{rmk:eta_p}
The function $\eta_p(x) \equiv \frac{p}{1-p}\eta(x^{1/p})x^{1-1/p}$ defined in \eqref{eq:defEta_p} appears (upto the factor $1-p$) naturally in connection with the deterministic Bihari-LaSalle equality: Let $x$ be as in \cref{lemma:DetBihari} and assume that $x$ is non-decreasing. Then, $x^p$ (for some $p\in(0,1)$) satisfies:
\begin{equation}\label{eq:4-detBihari-p}
x(t)^p \leq (1-p)\int_0^t\eta_p(x(s^-)^p)\dd A(s) + H^p  \quad \text{for all } t\in[0,T].
\end{equation}
This can be seen using the same formulas we will use in the stochastic case, see \cref{rmk:formulaZ}. \\

Moreover, for any continuous $\eta$ such that $x(t) \define G^{-1}(G(H)+t)$ is well-defined for $t\in[0,T]$, 
we have that $x$ satisfies
\begin{equation*}
x(t) = \int_0^t \eta(x(s)) \dd s + H
\end{equation*}
and $x^p$ satisfies
\begin{equation*}
x^p(t) = (1-p)\int_0^t \eta_p( x^p(s)) \dd s + H^p.
\end{equation*}
\end{remark}

\begin{remark}[On the relation between the concavity of $\eta_p$ and $\eta$] 
Assume that $\eta_p(0)=0$. Then, concavity of $\eta_p$ implies that $\eta$ is concave: 
If $\eta_p$ is concave, then it is almost everywhere differentiable. For any $x>0$ in which $\eta_p$ is differentiable, we have for $y=x^{1/p}$
\begin{equation*}
\eta(y)  =  \frac{1-p}{p}\eta_p(y^p) y^{1-p} \text{ and } \eta'(y) =\frac{1-p}{p}\left( p  \eta_p'(y^p) + (1-p)\eta_p(y^p)y^{-p}\right).
\end{equation*}
Due to $\eta_p$ being concave and $\eta_p(0) = 0$, we have that $z\mapsto\frac{\eta_p(z)}{z}$
is non-increasing.
In particular, $\eta$ is  almost everywhere differentiable and $\eta'$ is non-increasing (on its domain).  As concavity of $\eta_p$ implies, that $\eta_p$ (and hence also $\eta$) is locally Lipschitz continuous, we have that $\eta$ is absolutely continuous. Together, this implies that $\eta$ is concave. \\

However, concavity of $\eta$ does not imply that $\eta_p$ is concave: For example $\eta(x) = x \arctan(1/x)$ is concave (due to $\eta''(x) = -2(x^2+1)^{-2} <0$) and $\eta_{1/2}(x) = x \arctan(1/x^2)$ is not concave (as $\eta_{1/2}''(x) = 2x(x^4-3)(x^4+1)^{-2}$).
\end{remark}

\begin{remark}[Well-definedness of the upper bounds] 
In the deterministic Bihari-LaSalle inequality \cref{lemma:DetBihari} the assumption $G(H)+A(t) \in \textrm{domain}(G^{-1})$ is needed. In \cref{thm:stochBihariConcave} we do not need a corresponding assumption, because this is automatically satisfied if $\eta$ or $\eta_p$ is concave: We have $\dom = \text{range}(G)= (\lim_{\varepsilon \to 0} G(\varepsilon), \lim_{x\to\infty} G(x))$. Hence, it suffices to show $\lim_{x\to\infty}G(x) = \infty$. Concavity of $\eta$ implies that there exists some $K>0$ such that $\eta(u) \leq K u$ for all $u\in[c,\infty)$, and hence 
\begin{equation*}
\lim_{x\to\infty} G(x) = \int_c^{\infty} \frac{\dd u}{\eta(u)} \geq   \int_c^\infty  \frac{\dd u}{ K u} = +\infty.
\end{equation*}
Due to \eqref{eq:G_umrechnen} a similar argument implies $\lim_{x\to\infty}G(x) = \infty$ if $\eta_p$ is concave.
\end{remark}

\begin{remark}[Comparison of bounds in the cases $\Anosup$ and $\Asup$]
Recall that Assumption $\Anosup$ is a special case of Assumption $\Asup$. The bound given by 
\cref{thm:stochBihariConcave} b) is indeed better than the bound of \cref{thm:stochBihariConcave} a). This can be seen as follows: Under the assumptions of  \cref{thm:stochBihariConcave} a) we show for all $h>0, x\geq 0$ 
\begin{equation}
\alpha_1 G^{-1}(G(h)+x) \leq 
G^{-1}(G(\alpha_1 h) + \beta x)
\end{equation}
holds true. Fix some $h>0$. Define for all $x\geq 0$:
\begin{equation*}
\begin{aligned}
f(x) & \define  \alpha_1^p \big(G^{-1}(G( h)+ x)\big)^p = \alpha_1^p (\tilde{G}_p^{-1}( \tilde{G}_p(h^p)  + (1-p) x)), \\
g(x) &\define \big(G^{-1}(G(\alpha_1 h) + \beta x)\big)^p
=  \tilde{G}_p^{-1}(\tilde{G}_p(\alpha_1^p h^p) + (1-p)\beta x).
\end{aligned}
\end{equation*}
Here, we used \eqref{eq:G_umrechnen} and \eqref{eq:G_umrechnen2}. Note that $f(0) = \alpha_1^p h^p = g(0)$. Furthermore, we have $\alpha_1^p = \beta = (1-p)^{-1}$, $\alpha_1^{-p} \in(0,1)$ and
\begin{equation*}
\begin{aligned}
f'(x) & = (1-p)\alpha_1^p\frac{1}{\tilde{G}_p'(\alpha_1^{-p}f(x))}
=  \eta_p(\alpha_1^{-p}f(x)) \\
g'(x)& = (1-p)\beta \frac{1}{\tilde{G}_p'(g(x))} =  \eta_p(g(x)).
\end{aligned}
\end{equation*}
For concave $\eta_p$ the ODEs $y' = \eta_p(\alpha_1^{-p} y)$ and $y' = \eta_p(y)$ have unique solutions for $y(0)>0$. Therefore, the assumption that $\eta_p$ is non-decreasing implies that $f(x) \leq g(x)$ for all $x\geq 0$. Taking the $p$-th root implies the claim.
\end{remark}

%%%%%%%%%%%%%%%%%%%%%%%%%%%%%%%%%%%%%%%%%%%%%%%%%%%%%%%%%%%%%%%%%%%%%%%%%%%%%%%%%%%%%%%%%%

\subsection{Estimates for concave $\eta$ and random $A$}
In this section we extend the results of the previous subsection to predictable integrators $A$. Note that we may not expect in general estimates of the type
\begin{equation*}
\|X^*_T\|_p  \leq c_1 \big\|G^{-1}\big(G\big(c_2 H_T \big) + A_T\big)\big\|_p
\end{equation*}
to hold, where $c_1$ and $c_2$ are constants that only depend on $p\in(0,1)$. This type of estimate fails to hold true even for constant $H$, $\eta(x) \equiv x$ and  assumption $\Anosup$, see \cref{ex:counterexample2} in the appendix.
Due to the estimates known for the linear case $\eta(x) \equiv x$ (see e.g. \cite[Corollary 5.4.]{GeissConvex} and the references therein), we conjecture that estimates of the form
\begin{equation}\label{eq:WunschAbschaetzung}
\| X^*_t \big(G^{-1}(G(c_3 H)+A_T)\big)^{-1}\|_p \leq c_4,
\end{equation}
(where $G^{-1}$ denotes the inverse of $G$ and $(...)^{-1}$ the reciprocal) for some constants $c_3$ and $c_4$ that only depend on $p\in(0,1)$, might hold true for constant $H_t \equiv H$ and \nameref{def:nosup} (and similarly for \nameref{def:sup}).

Using \cref{thm:stochBihariConcave} we prove a slightly weaker estimate which is presumably not sharp.

\begin{corollary}[A stochastic Bihari-LaSalle inequality for concave $\eta$ and random $A$]\label{cor:randomA} 
Let $(X, A, H, M)$  and $\eta$ satisfy Assumption $\Asup$ (see \cref{def:sup}) assume that $\eta_p(x) \define \tfrac{p}{1-p} \eta(x^{1/p})x^{1-1/p}$ is concave and non-decreasing. Let $t>0$.
\begin{enumerate}
\item If $\E[H^p\ti] < \infty$ and either $H$ is predictable or $M$ has no negative jumps, then 

\begin{equation*}
\EO[\big(G^{-1}(G(\alpha_1 \alpha_2 \|H_t\|\pF) + \beta A_t)\big)^{-q} (X_t^*)^p]  \leq \frac{q}{q-p} (\alpha_1 \alpha_2 \|H_t\|\pF)^{p-q}
\end{equation*}
holds for all $p\in(0,1)$, and $q>p$. In particular, by an application of  H\"older's inequality  we also have
\begin{equation*}
\|X^*_t\|_{p/2,\F_0} \leq \left(\frac{q}{q-p} (\alpha_1 \alpha_2 \|H_t\|\pF)^{p-q}\right)^{1/p} \EO[\big(G^{-1}(G(\alpha_1 \alpha_2 \|H_t\|\pF) + \beta A_t)\big)^{q}]^{1/p}.
\end{equation*}

\item If $\E[H\ti] < \infty$ then 
\begin{equation*}
\EO[\big(G^{-1}(G(\alpha_1 \|H_t\|_{1,\F_0}) + \beta A_t)\big)^{-q} (X_t^*)^p]  \leq \frac{q}{q-p} (\alpha_1 \|H_t\|_{1,\F_0})^{p-q}
\end{equation*}
holds for all $p\in(0,1)$, and $q>p$. In particular, we also have
\begin{equation*}
\|X^*_t\|_{p/2,\F_0} \leq \left(
\frac{q}{q-p} (\alpha_1 \|H_t\|_{1,\F_0})^{p-q} \right)^{1/p} \EO[\big(G^{-1}(G(\alpha_1 \|H_t\|_{1,\F_0}) + \beta A_t)\big)^{q}]^{1/p}.
\end{equation*}
\end{enumerate}
If we assume the stronger assumption $\Anosup$ then we may replace $\beta$ by $1$ in the estimates above.
\end{corollary}

Moreover, under the stronger assumption $\Anosup$ we obtain:
\begin{corollary}[A stochastic Bihari-LaSalle inequality for concave $\eta$ and random $A$ for $\Anosup$ ]\label{cor:randomAnoSup}
Let $(X, A, H, M)$ satisfy Assumption $\Anosup$ (see \cref{def:nosup}) and assume that $\eta$ is concave. Then, for all $t\geq 0$, $R>0$, $w> 0$, $u> 0$ we have
\begin{equation*}
\PO[X^*_t > u] \leq 
\begin{cases}
\frac{G^{-1}(G(\EO[H_t\wedge w]) + R)}{u} + \PO[H_t \geq w] + \PO[A_t >  R]
 &\text{ if } H \text{ is predictable and } \eta \text{ is concave} \\[1em]
\frac{G^{-1}(G(\EO[H_t]) + R)}{u} + \PO[A_t >  R]
 &\text{ if } \E[H\ti] < \infty \text{ and } \eta \text{ is concave}.
\end{cases}
\end{equation*}
\end{corollary}

\cref{cor:randomAnoSup} does not hold under the weaker assumption $\Asup$, see \cref{rmk:noWeak}.

\subsection{Estimates for general $\eta$ and random $A$}
Now we study upper bounds of $X$ without imposing concavity or convexity assumptions on $\eta$. Whilst we do not obtain bounds for $\|X_T^*\|_p$ in this case, we do observe the same qualitative behavior as in the deterministic case.

\begin{remark}[Qualitative behavior of the deterministic Bihari-LaSalle inequality]\label{rmk:qualitativeBehavior}
Let $x$, $H$, $A$, $\eta$ and $G$ be as in \cref{lemma:DetBihari}.
\begin{enumerate}
\item Assume $\int_0^{\varepsilon} \frac{\dd u}{\eta(u)}=+\infty$ for all $\varepsilon>0$  (Osgood condition) and assume that 
\eqref{eq:detBihari} already holds when $H$ is replaced by $0$. The assumption on $\eta$ implies $\lim_{\varepsilon \to 0} G(\varepsilon) = -\infty$, and hence for small enough $\varepsilon>0$, we have that $G(\varepsilon)+A(t) \in \textrm{domain}(G^{-1})$. Hence, for any $\varepsilon>0$ we may apply \cref{lemma:DetBihari} to $x$, $\varepsilon$ (instead of $H$), $A$ and $\eta$ and obtain $x(t) \leq  \lim_{\varepsilon\to 0} G^{-1}(G(\varepsilon)+A(t)) = 0$.
\item If $\int_1^{\infty} \frac{\dd u}{\eta(u)}=+\infty$, then $[G(H),\infty) \subseteq \textrm{range}(G) 
=\textrm{domain}(G^{-1})$. This implies that 
$G(H) + A(t)\in\textrm{domain}(G^{-1})$, so we may apply \cref{lemma:DetBihari} to obtain the upper bound $x(t) \leq G^{-1}(G(H)+A(t))$.
\end{enumerate}
\end{remark}
\cref{thm:BihariRandomA} provides a stochastic version of
\cref{rmk:qualitativeBehavior}. Assertion \ref{item:SB-4} can be seen as an extension of \cite[Lemma 5.1]{RenesseScheutzow} from continuous to \cadlagg local martingales $M$.

\begin{theorem}[Stochastic Bihari-LaSalle type estimates for $\Asup$]\label{thm:BihariRandomA}
Let $(X, A, H, M)$ and $\eta$ satisfy assumption $\Asup$ and assume $\lim_{x\to\infty}\eta(x) = +\infty$. Then the following assertions hold.
\begin{enumerate}
\item\label{item:SB-1} We have for all stopping times $\tau$ and all deterministic times $T\geq 0$
\begin{equation*}
\E[G(X_{\tau\wedge T})] \leq \E[A_{\tau\wedge T}] + G(\E[H_{\tau\wedge T}]),
\end{equation*}
where the terms take values in $[0,+\infty]$. For the definition of $G$ see \eqref{eq:defG}.
\item\label{item:SB-2}  If $\int_0^{\varepsilon} \frac{\dd u}{\eta(u)}=+\infty$ for all $\varepsilon>0$ and $H_T=0$ then $X_T^*=0$ $\PP$-almost surely. 
\item\label{item:SB-3}  Let $(X^{(n)}, A, H^{(n)}, M^{(n)})$ satisfy $\Asup$ for each $n\in\N$ for some $\eta$ which does not depend on $n$. If $\int_0^{\varepsilon} \frac{\dd u}{\eta(u)}=+\infty$ for all $\varepsilon >0$ and $\lim_{n\to\infty}\E[H_T^{(n)}]=0$, then 
\begin{equation*}
\lim_{n\to\infty} \PP[\sup_{t\in[0,T]}X^{(n)}_t > \delta] = 0 \quad \text{ for each } \delta >0.
\end{equation*}
\item\label{item:SB-4}  If $X_t\geq c$ for all $t\in[0,T]$ for some $T>0$ where $c$ is the constant from the definition of $G$, then for any $p\in(0,1)$
\begin{equation*}
\|G(X_T^*)\|_p \leq \alpha_1\alpha_2\|A_T+ G(\E[H_T])\|_p
\end{equation*} 
where $\alpha_1$, $\alpha_2$ are constants that only depend on $p$ (see \eqref{eq:constants}). If we assume in addition that  $\int_1^{\infty} \frac{\dd u}{\eta(u)}=+\infty$, $\E[A_T^p]<\infty$ and $\E[H_T]<\infty$ this guarantees that $X_T^*$ is almost surely finite.
\end{enumerate}
\end{theorem}

\cref{thm:BihariRandomA} can be applied to study existence and uniqueness of global solutions to path-dependent SDEs:
In \cref{sec:applications} we will use \cref{thm:BihariRandomA}\ref{item:SB-2} to prove the uniqueness of the solutions, \cref{thm:BihariRandomA}\ref{item:SB-3} will be applied to prove that the Euler approximates are indeed a Cauchy sequence. Moreover, \cref{thm:BihariRandomA}\ref{item:SB-4} is useful to prove the non-explosion of the solution.

%%%%%%%%%%%%%%%%%%%%%%%%%%%%%%%%%%%%%%%%%%%%%%%%%%%%%%%%%%%%%%%%%%%%%%%%%%%%%%%%%%%%%%%%%%
%%%%%%%%%%%%%%%%%%%%%%%%%%%%%%%%%%%%%%%%%%%%%%%%%%%%%%%%%%%%%%%%%%%
%%%%%%%%%%%%%%%%%%%%%%%%%%%%%%%%%%%%%%%%%%%%%%%%%%%%%%%%%%%%%%%%%%%
%%%%%%%%%%%%%%%%%%%%%%%%%%%%%%%%%%%%%%%%%%%%%%%%%%%%%%%%%%%%%%%%%%%

\section{Application: Existence and uniqueness of solutions to a \levy driven path-dependent SDE}
\label{sec:applications}

The first stochastic Gronwall inequality \cite{RenesseScheutzow} (of the class of stochastic Gronwall inequalities related to the results of this paper) was developed to study stochastic functional differential equations: Stochastic Gronwall inequalities seem necessary to obtain the existence and uniqueness of local and global solutions under a one-sided Lipschitz condition.
The results of \cite{RenesseScheutzow} were extended to path-dependent SDEs driven by a more general integrator in \cite{MehriScheutzow}. In this section, using the stochastic Bihari-LaSalle inequality \cref{thm:BihariRandomA}, we obtain  existence and uniqueness of global solutions to path-dependent SDEs under a one-sided \textit{non-Lipschitz} condition. \\

Assume an underlying filtered probability space $(\Omega, \F,\PP,(\F_t)_{t\geq 0})$ satisfying the usual conditions. Let $L$ be a $\R^d$-valued \cadlagg \levy process with bounded jumps and its L\'evy-It\^o decomposition given by 
\begin{equation*}
L\ti = bt + \sigma B\ti + \int_{|\xi| \leq c} \xi \tilde{N}(t,\dd \xi) \qquad  \forall t\geq 0
\end{equation*}
where $b\in\R^d$, $\sigma \in \R^{d\times m}$, $B$ an $m$-dimensional Brownian motion, $N$ an independent Poisson random measure  on $\R^+\times (\R^d-\{0\})$ with \levy measure $\nu$  and $c>0$ some constant. We denote by $\tilde{N}$ the  corresponding compensated Poisson random measure. For details, see for example \cite[p. 126]{Apfelbaum}. Let $| \cdot|_F$ denote the Frobenius norm on $\R^{d\times m}$. Define $U\define \{ \xi\in\R^d \mid 0 < | \xi | \leq c\}$ and set $\mathcal{U}\define \B(U)$, where $\B(U)$ denotes the Borel $\sigma$-algebra on $U$.

\bigskip

\noindent We study the path-dependent SDE with random coefficients driven by $L$:
\begin{equation}\label{eq:SDE}
\begin{cases}
d X_t & = f(t,X) \dd t +  \int_{|\xi| \leq c} g(t, X, \xi) \tilde{N}(\dd t, \dd \xi)  + h(t, X) \dd B\ti \\
X_t & = z_t, \qquad t\in[-r, 0],
\end{cases}
\end{equation}
where $r>0$ is some constant and the initial condition $(z_t)_{t\in[-r,0]}$  has \cadlagg paths and is $\F_0$ measurable. Denote by $\pred$ the predictable $\sigma$-algebra on $[0,\infty)\times\Omega$. 
Recall that the Borel $\sigma$-field on the space $\cadlag$ induced by the Skorohod metric $J_1$ is strictly contained in the Borel $\sigma$-field generated by the uniform norm, see e.g. \cite[Section 6, Setion 12]{Billingsley} for more details. Whilst we work on the space $\cadlag$ predominantly with uniform convergence, we need the coefficients $f$, $g$, $h$ to be measurable w.r.t. the Borel $\sigma$-field induced by convergence in Skorohod metric $J_1$ on compacts sets, which we denote by $\tilde{\B}(\cadlag))$. We need this stronger measurability assumption to ensure the measurability of the coefficients in \eqref{eq:SDE}.

Assume that the coefficients
\begin{equation*}
\begin{aligned}
& f:([0,\infty)\times  \Omega \times \cadlag, \pred \otimes \tilde{\B}(\cadlag)) \to (\R^d,\B(\R^d)), \\
& g:([0,\infty) \times \Omega \times \cadlag \times U, \pred \otimes \tilde{\B}(\cadlag))\otimes \mathcal{U} \to (\R^d,\B(\R^d)), \\
& h:( [0,\infty) \times \Omega \times\cadlag, \pred \otimes \tilde{\B}(\cadlag)) \to (\R^{d\times m},\B(\R^{d\times m})) \\
\end{aligned}
\end{equation*}
are measurable mappings. For every $t\in[0,\infty)$, $\omega\in\Omega$, $\xi\in U$ assume that $f(t, \omega, x)$, $g(t, \omega, x,\xi)$ and $h(t, \omega, x)$ only depend on the path segment $x(s),s\in[-r, t)$. We denote by $f(t,x)$, $g(t,x,\xi)$ and $h(t,x)$ the corresponding random variables. \\

\begin{hypo}\label{hypo}
Let  $(L^R_t)_{t\geq 0}$ and $(\tilde{K}^R_t)_{t\geq 0}$ for all $R>0$ and $(K_t)_{t\geq 0}$ be non-negative, jointly measurable and adapted stochastic processes, such that for all $t\geq 0$ the quantities $\E[\int_0^t L^R_s \dd s]$, $\E[\int_0^t \tilde{K}^R_s \dd s]$ and $\E[\int_0^t K_s \dd s]$  are finite. Assume that for all $\omega\in\Omega$, for all $x,y\in\text{C\`adl\`ag}([-r, \infty);\R^d)$ and for all $t\geq 0$:
%There exist non-negative functions $K\in L_{loc}^1([0,\infty),\dd t)$ and $L_R, 
%\tilde{K}_R \in L_{loc}^1([0,\infty),\dd t)$ for all $R>0$ 
\begin{enumerate}
\item[(C1)]\label{item:c1} for $\sup_{s\in[-r, t]} |x(s)|$,  $\sup_{s\in[-r, t]} |y(s)|  \leq R$
\begin{equation*}
\begin{aligned}
2 \langle x(t^-) - y(t^-), f(t,\omega, x) - f(t, \omega, y) \rangle + \int_{U} |g(t,\omega, x, \xi) - g(t,\omega, y,\xi)|^2 \nu (\dd \xi)  \\ 
+ |h(t,\omega, x) - h(t,\omega, y)|_F^2 \leq L^R_t(\omega) \eta_1( \sup_{s\in[-r,t]}|x(s)-y(s)|^2),
\end{aligned}
\end{equation*}

\item[(C2)\label{item:c2}]  \begin{minipage}[t][.5cm][b]{0.9\textwidth}
\begin{equation*}
\begin{aligned}
 2 \langle x(t^-), f(t,\omega, x) \rangle + \int_{U} |g(t,\omega, x, \xi)|^2 \nu (\dd \xi) 
 + |h(t,\omega, x)|_F^2  \\ \leq K_t(\omega) \eta_2(1 + \sup_{s\in[-r,t]}|x(s)|^2),
 \end{aligned}
 \end{equation*}
\end{minipage} 
%\begin{equation*}
%\begin{aligned}
% 2 \langle x(t^-), f(t,\omega, x) \rangle + \int_{U} |g(t,\omega, x, \xi)|^2 \nu (\dd \xi) 
% + |h(t,\omega, x)|_F^2  \\ \leq K_t(\omega) \eta_2(1 + \sup_{s\in[-r,t]}|x(s)|^2),
% \end{aligned}
% \end{equation*}
%\item[(C2)]\label{item:c2} \vspace{-2.2cm}
%
%\phantom{aaaa}\bigskip\bigskip\bigskip\bigskip

\item[(C3)]\label{item:c3} $x\mapsto f(t,\omega, x)$ as a function from $\text{C\`adl\`ag}([-r, \infty);\R^d)$ (endowed with the topology induced by uniform convergence on compact sets) to $\R^d$ is continuous,
\item[(C4)]\label{item:c4} for $\sup_{s\in[-r, t]} |x(s)| \leq R$
\begin{equation*}
|f(t,\omega, x)| + \int_{U} |g(t,\omega, x, \xi)|^2 \nu (\dd \xi) + |h(t,\omega, x)|_F^2 \leq \tilde{K}^R_t(\omega),
\end{equation*}
\item[(C5)]\label{item:c5} $\E[\sup_{s\in[-r,0]} |z_s|^2] < \infty$, 
\end{enumerate}
where $\eta_1,\eta_2:[0,\infty) \to [0,\infty)$ are non-decreasing continuous functions with $\eta_i(x) >0$ for all $x>0$, $i=1,2$. Define 
\begin{equation*}
G_{i}(x) \define \int_1^x \frac{\dd u}{\eta_{i}(u)}, \qquad \forall x>0, i=1,2.
\end{equation*}
Assume:
\begin{equation*}
\begin{aligned}
\int_0^\varepsilon\frac{\dd u}{\eta_{1}(u)} & =+\infty   \quad \forall \varepsilon>0  & \text{ i.e. } \lim_{\varepsilon\searrow 0} G_1(\varepsilon) = -\infty, \\
\int_1^\infty \frac{\dd u}{\eta_2(u)} &=+\infty &\text{ i.e. }  \lim_{x\to\infty} G_2(x)=+\infty.
\end{aligned}
\end{equation*}
\end{hypo}

\begin{remark}
The assumption $\lim_{\varepsilon\searrow 0} G_1(\varepsilon) = -\infty$ implies that $\lim_{\varepsilon\searrow 0} \eta_1(\varepsilon) = 0$. Hence, (C1) and (C3) imply that for any fixed $\omega\in\Omega$ and $t\geq 0$ we have that the mappings
\begin{equation*}
\begin{aligned}
&\cadlag \to L^2(U, \mathcal{U}, \nu), \qquad &&x\mapsto g(t,\omega,x, \xi),  \\
&\cadlag \to \R^{d\times m},  \qquad&& x\mapsto h(t,\omega, x) 
\end{aligned}
\end{equation*}
are continuous  (with respect to topology induced by 
uniform convergence on compact sets).
\end{remark}

The following theorem is a nonlinear extension of \cite[Theorem 3.3]{MehriScheutzow}:
We study an SDE driven by a less general integrator than \cite{MehriScheutzow}, but we have weaker assumptions on the coefficients, as we allow in (C1) and (C2) of \cref{hypo} nonlinear upper bounds $\eta_1$ and $\eta_2$. We also impose slightly weaker integrability assumptions on the coefficients. Our proof of existence and uniqueness of global solutions is similar to that of \cite{MehriScheutzow} except that we use a stochastic Bihari-LaSalle inequality (\cref{thm:BihariRandomA}) instead of a stochastic Gronwall inequality. In contrast to \cite{RenesseScheutzow} our SDE has \cadlagg paths and we have a nonlinear monotonicity condition (C1). 

\begin{corollary}[Existence and uniqueness of global of solutions]\label{cor:solutions}
Assume \cref{hypo} holds. Then the SDE \eqref{eq:SDE} has a unique strong global solution.
\end{corollary}

The following example satisfies the assumptions of \cref{cor:solutions} but does not seem to be 
covered by the results in the existing literature. It is similar to \cite[Example 4.1]{MehriScheutzow}.
\begin{example}
Consider the following $1$-dimensional path-dependent (or functional) SDE
\begin{equation*}
\dd X_t = \bigg[ -2\textrm{sgn}(X_{t^-})|X_{t^-}|^{1/2} + \sup_{t-1\leq s<t}|X_s|\bigg] 
\dd t + \bigg[ |X_{t^-}|^{3/4} +  \sup_{t-1\leq s<t}|X_s| + \big(|X_{t^-}|\wedge \tfrac{1}{\e}\big) \log\bigg(\tfrac{1}{|X_{t^-}|\wedge \tfrac{1}{\e}}\bigg)\bigg] \dd W_t + \dd N_t.
\end{equation*}
where $B$ is a standard Brownian motion and $N$ a compensated standard Poisson process.

Due to the term  $-\big(|X_{t^-}|\wedge \tfrac{1}{\e}\big) \log\big(|X_{t^-}|\wedge \tfrac{1}{\e}\big)$ the diffusion coefficient does not satisfy a one-sided Lipschitz assumption, hence the results of \cite{RenesseScheutzow} and \cite{MehriScheutzow} cannot be applied. Moreover, due to the path-dependence or the jump term $\dd N_t$ the results of \cite{WuLan} cannot be applied.

As mentioned in the introduction, existence and uniqueness results could be proven without stochastic Bihari-LaSalle inequalities (using the BDG inequality instead) if the inequality of (C1) is replaced by the stronger assumption
\begin{equation*}
\begin{aligned}
2 \langle x(t^-) - y(t^-), f(t,\omega, x) - f(t, \omega, y) \rangle \leq L^R_t(\omega) \eta_1( \sup_{s\in[-r,t]}|x(s)-y(s)|^2) \text{ and }\\
\int_{U} |g(t,\omega, x, \xi) - g(t,\omega, y,\xi)|^2 \nu (\dd \xi) 
+ |h(t,\omega, x) - h(t,\omega, y)|_F^2 \leq L^R_t(\omega) \eta_1( \sup_{s\in[-r,t]}|x(s)-y(s)|^2).
\end{aligned}
\end{equation*}
Due to the term $|X_{t^-}|^{3/4}$ in the diffusion coefficient, this strengthened assumption is not satisfied by the diffusion coefficient $g$ for suitable $\eta_1$.

It is easily seen that the coefficients satisfy assumptions (C2), (C3) and (C4) of \cref{hypo}. 
We show that also (C1) is satisfied:
We have for $0<x<y$:
\begin{equation*}
|y^{3/4} - x^{3/4}|^2 = (3/4)^2\big(\int_x^y s^{-1/4} \dd x\big)^2 \leq (y-x) \int_x^y s^{-2/4} \dd s  \leq  2 (y-x)(y^{1/2}-x^{1/2}).
\end{equation*}
This implies for all $x,y \in \cadlag$ that
\begin{equation*}
2\big\langle y(t^-) - x(t^-),-2\textrm{sgn}(y(t^-))|y(t^-)|^{1/2} +2\textrm{sgn}(x(t^-))|x(t^-)|^{1/2} \big\rangle  + 2\big| |y(t^-)|^{3/4}- |x(t^-)|^{3/4}\big|^2 \leq 0.
\end{equation*}

Define $\varphi(x) \define y\log(1/y) -x\log(1/x) - (y-x)\log(1/(y-x))$ for $x\in[0,y]$. Noting that $\varphi$ is convex and $\varphi(0) = \varphi(y) = 0$, implies 
$g\leq 0$ and in particular for all $0\leq x\leq y\leq 1/\e$
\begin{equation*}
|y\log(1/y) - x\log(1/x)| \leq |y-x|\log(1/|y-x|).
\end{equation*}
Hence, choosing $\eta_1(x) = K(x+(x \wedge \e^{-1})\log\big(\frac{1}{x \wedge \e^{-1}}\big))$ for some $K>0$ large enough implies that the coefficients indeed satisfy (C1).
\end{example}

%%%%%%%%%%%%%%%%%%%%%%%%%%%%%%%%%%%%%%%%%%%%%%%%%%%%%%%%%%%%%%%%%%%%
%%%%%%%%%%%%%%%%%%%%%%%%%%%%%%%%%%%%%%%%%%%%%%%%%%%%%%%%%%%%%%%%%%%%

\section{Proofs}

\subsection{A Lenglart type estimate}

We will prove the stochastic Bihari-LaSalle inequalities using the following lemma which is \cite[Lemma 4.5]{GeissConvex}. It is an extension of an inequality Lenglart inequality, see \cite[Th\'eor\`eme I]{Lenglart} or \cite[Lemma 2.2 (ii)]{MehriScheutzow}:

\begin{lemma}[Lenglart type estimates]\label{lemma:1}
Fix some $T>0$ and $p\in(0,1)$ and let $(X,A,H,M)$ satisfy \nameref{def:sup} or \nameref{def:nosup}. We consider the following 6 cases, which arise from combining $\Asup$ or $\Anosup$ with one of the following three assumptions:
\begin{enumerate}
\item $H$ is predictable and $\E[H^p_T]<\infty$,
\item $M$ has no negative jumps and $\E[H^p_T]<\infty$,
\item $\E[H_T]<\infty$.
\end{enumerate}
Fix arbitrary $u,\lambda >0$  and set:
\begin{equation*}
\tau_u \define \tau \define \inf\{s \geq 0 \mid H\s \geq \lambda u\}, \qquad \sigma_u \define \sigma \define \inf\{ s\geq 0 \mid X\s > u\},
\end{equation*}
where $\inf\emptyset \define +\infty$. Then, the following estimate holds true for all $t\in[0,T]$:
\begin{equation}\label{eq:lemma1}
\one_{\{X^*_t > u\}} u \leq  X_{t\wedge \sigma_u} \wedge u  \leq I^{L,u}\ti + M^{L,u}\ti + H^{L,u}\ti.
\end{equation}
Here
$(I^{L,u}\ti)_{t\geq 0}$ is a non-decreasing process containing the integral term from \eqref{eq:GronwallAssumption} and \eqref{eq:GronwallAssumptionNoSup} respectively with an additional indicator function
\begin{equation*}
I^L\ti \define I^{L,u}\ti \define
\begin{cases}
\int_{(0,t]} \eta(X^*\sm) \one_{\{X^*\sm \, \leq u \}} \dd A\s  & \text{for } \Asup, \\[1em]
\int_{(0,t]} \eta(X\sm) \one_{\{X^*\sm \, \leq u \}} \dd A\s  & \text{for } \Anosup, 
\end{cases}
\end{equation*}
the process $(M^{L,u}\ti)_{t\geq 0}$ is a local martingale with \cadlagg paths starting in $0$ defined by
\begin{equation*}
M^L\ti \define M^{L,u}\ti \define
\begin{cases}
 \lim_{n\to\infty} M_{t\wedge \tau^{(n)} \wedge \sigma}  & \text{ if $H$ is predictable and } \E[H^p_T]<\infty,  \\
M_{t\wedge \tau \wedge \sigma} &  \text{ if $M$ has no negative jumps and } \E[H^p_T]<\infty, \\  
\tilde M_{t\wedge \sigma}\one_{\{\EO[H_T]\leq u\}}  & \text{ if } \E[H_T]<\infty,
\end{cases}
\end{equation*}
(where $\tau^{(n)}$ denotes an announcing sequence of $\tau$ and $\tilde M\ti \define M\ti + \E[H_T\mid F_t] -\EO[H_T]$ for $t\in[0,T]$), 
and $(H^{L,u}\ti)_{t\geq 0}$ is a non-decreasing process depending on $H$:
\begin{equation*}
H^L\ti \define H^{L,u}\ti \define 
\begin{cases}
H\ti \wedge (\lambda u)  + u\one_{\{H\ti \geq  \lambda u\}} & \text{ if $H$ is predictable and } \E[H^p_T]<\infty,  \\
H\ti \wedge (\lambda u)  + u\one_{\{H\ti \geq  \lambda u\}}  &  \text{ if $M$ has no negative jumps and } \E[H^p_T]<\infty, \\  
\EO[H_T]\wedge u &  \text{ if } \E[H_T]<\infty.  
\end{cases}
\end{equation*}
\end{lemma}

We will frequently use the formulas of the subsequent remark.
\begin{remark}[Calculation of $Z^p$, $p\in(0,1)$]\label{rmk:formulaZ}
Let $Z$ be a non-negative random variable and $p\in(0,1)$.
Then $Z^p$ can be calculated using the three formulas below.
\begin{equation}\label{eq:formulaZ}
\begin{aligned}
Z^p &= p \int_0^\infty \one_{\{Z \geq u\}} u^{p-1} \dd u \\ 
Z^p & = (1-p) \int_0^\infty Z \one_{\{Z\leq u\}} u^{p-2} \dd u \\
Z^p &= p (1-p) \int_0^\infty (Z\wedge u)u^{p-2}\dd u
\end{aligned}
\end{equation}
The third equality follows e.g. by using the first and second equality.
In particular, we also have  $Z^{p-1}  = (1-p) \int_0^\infty \one_{\{Z\leq u\}} u^{p-2} \dd u$ for $Z>0$.
\end{remark}

The third equality of \eqref{eq:formulaZ} holds true for general concave functions, which is due to Burkholder \cite[Theorem 20.1, p.38-39]{Burkholder}, see also Pratelli \cite[p. 403]{Pratelli}:

\begin{remark}\label{rmk:Burkholder}[Calculation of $F(Z)$ for concave $F$]
Let $F\colon \R \mapsto \R$ be a concave function and let $F'$ its left-hand derivative. Then, for any non-negative random variable $Z$ we have
\begin{equation*}
F(Z) = F(0) + Z F'(\infty) - \int_{[0,\infty)} (Z\wedge u) \dd F'(u).
\end{equation*}
For the convenience of the reader we provide the proof: Recall that  left-hand derivative $F'$ of a concave function is non-increasing and left-continuous. By integration by parts (see e.g. \cite[Theorem 21.67 (iv), p.  419]{Hewitt}) we have
\begin{equation*}
\begin{aligned}
\int_{[0,\infty)} (Z\wedge u) \dd F'(u)
 &= \int_{[0,Z]} u \dd F'(u) + Z \int_{[Z,\infty)} \dd F'(u) \\
 &=   Z F'(Z+) - \int_{[0,Z]} F'(u) \dd u +  Z F'(\infty) - Z F'(Z+) \\
& = - F(Z) + F(0) + ZF'(\infty).
\end{aligned}
\end{equation*}
The last equality can be seen by e.g. using that concave functions are absolutely continuous and almost everywhere differentiable.
\end{remark}

%%%%%%%%%%%%%%%%%%%%%%%%%%%%%%%%%%%%%%%%%%%%%%%%%%%%%%%%%%%%%%%%%%%
%%%%%%%%%%%%%%%%%%%%%%%%%%%%%%%%%%%%%%%%%%%%%%%%%%%%%%%%%%%%%%%%%%%
%%%%%%%%%%%%%%%%%%%%%%%%%%%%%%%%%%%%%%%%%%%%%%%%%%%%%%%%%%%%%%%%%%%

\begin{remark}\label{rmk:HgeqEps} 
In the proofs of the stochastic Bihari-LaSalle inequalities, we will assume $X\geq \varepsilon$ and $H\geq \varepsilon$ for some $\varepsilon>0$ instead of $X\geq 0$ and $H\geq 0$.  We do this to ensure that terms like $G(X_t)$ or $G(H_t)$ are well-defined and to ensure that we may apply the deterministic Bihari-LaSalle inequality, see \cref{lemma:DetBihari}.
We may do this without loss of generality because we can add an arbitrary $\varepsilon>0$ to \eqref{eq:GronwallAssumption} (or similarly \eqref{eq:GronwallAssumptionNoSup}) and slightly weaken \eqref{eq:GronwallAssumption} (using that $\eta$ is non-decreasing) to obtain for all $t\in[0,T]$:
\begin{equation*}
(X\ti + \varepsilon) \leq \int_{(0,t]} \eta(X^*\sm+\varepsilon) \dd A\s + M\ti + (H\ti + \varepsilon) \qquad \PP\text{-a.s}.
\end{equation*}
Proving the assertions of the theorems for the processes $(X\ti + \varepsilon)_{t\geq 0}$, $(A\ti)_{t\geq 0}$, $(M\ti)_{t\geq 0}$ and $(H\ti + \varepsilon)_{t\geq 0}$ and then taking the limit $\varepsilon\to 0$ will imply the assertions for the general case $X\geq 0$, $H\geq 0$. 
\end{remark}

\subsection{Proofs for concave $\eta$ (\cref{thm:stochBihariConcave})}

\begin{proof}[Proof of \cref{thm:stochBihariConcave}, \ref{item:concaveSup}]
The sharpness of the constants follows from \cite[Theorem 3.5 and Theorem 3.7]{GeissConvex}. Assume w.l.o.g. that $X\geq \varepsilon$ and $H\geq \varepsilon$ for some $\varepsilon>0$ (see \cref{rmk:HgeqEps}).
Furthermore, we assume w.l.o.g. that $M$ is a martingale. This implies, that $(M^{L,u}_t)_{t\in[0,T]}$ from \cref{lemma:1} is  a martingale for any $u>0$.

Fix some $t\in[0,T]$. \cref{lemma:1} implies (under the assumption $\Asup$) for all $u>0$, $\lambda>0$:
\begin{equation*}
\one_{\{X^*\ti > u\}} u\leq \int_{(0,t]} \eta(X^*\sm)\one_{\{X^*\sm \leq u\}}\dd A\s 
+ M^{L,u}\ti  + H^{L,u}\ti
\end{equation*} 
using the notation $M^{L,u}$ and $H^{L,u}$ from \cref{lemma:1}. Fix some $\hat{u}>0$. We multiply the equation above with $pu^{p-2}\one_{\{u < \hat{u}\}}$, take the conditional expectation given $\F_0$ and integrate over $u$. This implies (using Fubini and $\{X_t^*>u\} \cap \{\hat{u} >u\} = \{X^*\ti\wedge\hat{u} > u\} $):
\begin{equation}\label{eq:proof-concave-1}
\begin{aligned}
\int_0^\infty \PO[X^*\ti\wedge\hat{u} > u] pu^{p-1}\dd u &
\leq \EO\bigg[  \int_{(0,t]} \eta(X^*\sm) \int_0^{\hat{u}} \one_{\{X^*\sm \leq u\}} p u^{p-2}  \dd u    \dd A\s \bigg]   \\
& \qquad + \int_0^{\hat{u}} \EO[H^{L,u}\ti] p u^{p-2} \dd u.
\end{aligned}
\end{equation}
Hence, by \cref{rmk:formulaZ} inequality \eqref{eq:proof-concave-1} can be slightly weakened to
\begin{equation}\label{eq:proof-concave-1b}
\begin{aligned}
\EO[ (X^*\ti\wedge \hat{u})^p] &
 = \int_0^\infty \PO[X^*\ti\wedge\hat{u} > u] pu^{p-1}\dd u \\
 & \leq \EO\bigg[  \int_{(0,t]} \eta(X^*\sm\wedge \hat{u})
 \frac{p}{1-p} \big((X^*\sm\wedge\hat{u})^{p-1} - \hat{u}^{p-1}\big) \dd A\s \bigg]   \\
& \qquad + \int_0^{\hat{u}} \EO[H^{L,u}\ti] p u^{p-2} \dd u \\
& \leq \frac{p}{1-p} \EO\bigg[  \int_{(0,t]} \eta(X^*\sm \wedge\hat{u}) (X^*\sm \wedge \hat{u})^{p-1} \dd A\s \bigg] + \int_0^\infty\EO[H^{L,u}\ti] p u^{p-2} \dd u.
\end{aligned}
\end{equation}
If $\E[H_T^p] <\infty$ and either $H$ is predictable or $\Delta M \geq 0$, we have (choosing $\lambda = p$ and applying \cref{rmk:formulaZ} and recalling $\alpha_1 = (1-p)^{-1/p} $, $\alpha_2 =p^{-1}$):
\begin{equation}\label{eq:Hp-calc}
\begin{aligned}
\int_0^\infty H^{L,u}\ti p u^{p-2} \dd u 
& = \int_0^\infty \big(H\ti \wedge (\lambda u)  + u\one_{\{H\ti \lambda^{-1}\geq u\}} \big) p u^{p-2} \dd u \\
& =  \lambda \int_0^\infty \big((H\ti\lambda^{-1}) \wedge u\big) p u^{p-2} \dd u
+ \int_0^\infty\one_{\{H\ti \lambda^{-1}\geq u\}}  p u^{p-1} \dd u \\
& = \lambda (1-p)^{-1}  (H_t \lambda^{-1})^p + (H_t \lambda^{-1})^p  \\
& = (p(1-p)^{-1} + 1)p^{-p} H_t^p \\
& = \alpha_1^p \alpha_2^p  H_t^p. 
\end{aligned}
\end{equation}
If $\E[H_T] <\infty$, we have: 
\begin{equation}\label{eq:Hp-calc2}
\int_0^\infty H^{L,u}\ti p u^{p-2} \dd u  = \int_0^\infty\big(\EO[H_T]\wedge u\big) p u^{p-2} \dd u =  (1-p)^{-1}  \EO[H_T]^p = \alpha_1^p  \EO[H_T]^p. 
\end{equation}
Furthermore, by assumption, $\eta_p(x)\equiv \frac{p}{1-p}\eta(x^{1/p}) x^{1-1/p}$ is concave and $A$ is deterministic. Therefore, \eqref{eq:proof-concave-1b} implies by Jensen's inequality for all $t\in[0,T]$:
\begin{equation}\label{eq:Lemma1-integrated}
\begin{aligned}
\EO[ (X^*\ti\wedge \hat{u})^p] 
  & \leq   \EO\bigg[  \int_{(0,t]} \eta_p\big((X^*\sm \wedge\hat{u})^p\big) \dd A\s \bigg] + \int_0^\infty\EO[H^{L,u}\ti] p u^{p-2} \dd u \\
& \leq \int_{(0,t]} \eta_p(\EO[(X^*\sm\wedge\hat{u})^p]) \dd A\s 
\\ & \quad + 
\begin{cases}
\alpha_1^p \alpha_2^p  \EO[H_T^p] & \text{ if } H \text{ is predictable or } \Delta M \geq 0, \\
\alpha_1^p  \EO[H_T]^p & \text{ if } \E[H_T] < \infty.
\end{cases}
\end{aligned}
\end{equation}
Since $H\geq \varepsilon$,  $\eta_p$ is non-decreasing and $s\mapsto \EO[(X^*\s\wedge \hat{u})^p]$ is c\`adl\`ag and takes values in $[0,\infty)$, we may apply the deterministic Bihari-LaSalle inequality (see \cref{lemma:DetBihari}) to the previous inequality. 
Recall the following definition and its properties (see 
\eqref{eq:defTildeG} and \eqref{eq:G_umrechnen}):
\begin{equation*}
\tilde G_p(x) \define 
\int_{c^{p}}^{x} \frac{\dd u}{\eta_p(u)}, \quad \text{satisfying}  \,\, \tilde G_p(x) = (1-p)G(x^{1/p}), \quad \tilde{G}_p^{-1}(x) = \big(G^{-1}\big(\tfrac{x}{1-p} \big)\big)^p.
\end{equation*}
For e.g. the case that $H$ is predictable or $M$ has non-negative jumps the deterministic Bihari-LaSalle inequality implies for all $\hat{u}>0$
\begin{equation*}
\begin{aligned}
\|X^*_T\wedge\hat{u}\|\pF & \leq \bigg\{\tilde{G}_p^{-1}\bigg( \tilde{G}_p\big(\alpha_1^p \alpha_2^p \EO[H^p_T] \big)  + A_T\bigg)\bigg\}^{1/p} \\
& = \bigg\{\tilde{G}_p^{-1}\bigg( (1-p) G\bigg(\alpha_1 \alpha_2 \|H_T\|\pF \bigg)  + A_T\bigg)\bigg\}^{1/p} \\
& = G^{-1}\bigg(G\bigg(\alpha_1 \alpha_2\|H_T\|\pF\bigg)  + \beta A_T\bigg). \\
\end{aligned}
\end{equation*}
Taking the limit $\hat{u}\to\infty$ implies the claim. A nearly identical argument implies the claim for also for the case $\E[H_T]<\infty$.
\end{proof}

\begin{proof}[Proof of \cref{thm:stochBihariConcave}, \ref{item:concaveNoSup}] 
The sharpness of the constants follows from \cite[Theorem 3.7]{GeissConvex}. We assume w.l.o.g. that $X\geq \varepsilon$ and $H\geq \varepsilon$ for some $\varepsilon>0$ (see \cref{rmk:HgeqEps}) and that $M$ is a martingale. As before, we define the following stopping times for some fixed $u,\lambda>0$:
\begin{equation*}
\tau_u := \inf\{s \geq 0 \mid H\s \geq \lambda u\}, \qquad \sigma_u := \inf\{ s\geq 0 \mid X\s > u\}.
\end{equation*}
We first prove the assertion for the case $\E[H_T] <\infty$, which is related to the proof  of \cite[Lemma 3.7]{XieZhang}.\\
\textbf{Proof for $\E[H_T]<\infty$:}
As $\eta$ is concave and $H$ is non-decreasing, we have for all $t\in[0,T]$:
\begin{equation*}
\begin{aligned}
\EO[X_{t\wedge \sigma_u}\wedge u] & \leq \EO[\int_{(0,t\wedge\sigma_u]} \eta(X\sm)\dd A\s] +  \EO[H_T] \\
& \leq \EO[\int_{(0,t]} \eta(X_{s^-\wedge \sigma_u}\wedge u)\dd A\s] +  \EO[H_T] \\
& \leq \int_{(0,t]} \eta(\EO[X_{s^-\wedge \sigma_u}\wedge u])\dd A\s +  \EO[H_T]
\end{aligned}
\end{equation*}
Recall that $X$ is a non-negative process. Therefore, dominated convergence implies that $t\mapsto \E[X(t\wedge\sigma_u)\wedge u]$  is c\`adl\`ag. Thus, we may apply the deterministic Bihari-LaSalle inequality \cref{lemma:DetBihari}, which gives:
\begin{equation*}
\PO[X^*_T>u]u \leq \EO[X_{T\wedge \sigma_u}\wedge u]\leq G^{-1}(G(\EO[H_T]) + A_T)=:\delta.
\end{equation*}
This implies the assertion \eqref{eq:claimConcaveNoSup-weak}.
In particular, we have $\PO[X^*_T>u]u  \leq  \delta \wedge u$, and hence by \cref{rmk:formulaZ}
\begin{equation*}
\begin{aligned}
\|X^*_T\|^p_{p,\F_0} 
& = p \int_0^\infty \PP[X^*_T>u] u^{p-1} \dd u   
\leq p \int_0^\infty (\delta \wedge u)u^{p-2} \dd u  \\
 &  = \tfrac{1}{1-p} \delta^p = \tfrac{1}{1-p}  G^{-1}(G(\EO[H_T]) + A_T)^p,
\end{aligned}
\end{equation*}
which implies the claim \eqref{eq:claimConcaveNoSup-p} (recalling $\alpha_1 = (1-p)^{-1/p}$). \\ 

\noindent\textbf{Proof for predictable $H$ and $\Delta M \geq 0$:}
Assume that $\E[H_T^p] <\infty$ and that  either $H$ is predictable or $\Delta M\geq 0$.  As we assumed w.l.o.g. that $M$ is a martingale, also $(M^{L,u}_t)_{t\in[0,T]}$ from \cref{lemma:1} is  a martingale for any $u>0$. By \cref{lemma:1} (using the case $\Anosup$) we have for all $t\in[0,T]$
\begin{equation*}
\begin{aligned}
\EO[X_{t\wedge \sigma_u} \wedge u] & \leq 
\EO[\int_{(0,t]} \eta(X\sm) \one_{\{X^*\sm \leq u \}} \dd A\s ]\\ 
& + \lambda \EO[H\ti\lambda^{-1} \wedge u] + u \PO[H\ti\lambda^{-1} \geq u].
\end{aligned}
\end{equation*}
Integrating w.r.t to $p(1-p)u^{p-2} \dd u$, using \cref{rmk:formulaZ} and choosing $\lambda = p$ (see also \eqref{eq:Hp-calc}) gives:
\begin{equation*}
\begin{aligned}
y(t) & \define p(1-p)\int_0^\infty \EO[X_{t\wedge \sigma_u}\wedge u] u^{p-2}\dd u  \\
& \leq p \EO[\int_{(0,t]} \eta(X\sm) (X^*\sm)^{p-1} \dd A\s] + p^{-p} \EO[H^p\ti].
\end{aligned}
\end{equation*}
Using $p\in(0,1)$, we have $(X^*\sm)^{p-1} \leq X\sm^{p-1}$. This implies (using the definitions $\eta_p(x) = \frac{p}{1-p}\eta(x^{1/p})x^{1-1/p}$, $\alpha_2 = p^{-1}$) for all $t\in[0,T]$
\begin{equation*}
\begin{aligned}
y(t) & \leq (1-p) \int_{(0,t]} \EO[\eta_p(X\sm^p)]\dd A\s + \alpha_2^p\EO[H^p\ti].
\end{aligned}
\end{equation*}
Note, that we have by \cref{rmk:formulaZ} $\E[X_t^p] = p(1-p) \int_0^\infty \E[ X_t \wedge u] u^{p-2} \dd u \leq y(t)\leq \E[(X_t^*)^p]$. In particular, $y(t)<\infty$ for all $t\in[0,T]$ due to \cref{thm:stochBihariConcave}\ref{item:concaveSup}.  Hence, due to $\eta_p$ being concave and non-decreasing and $\E[X^p\ti] \leq y(t)$, we obtain:
\begin{equation*}
\begin{aligned}
y(t) & \leq  (1-p) \int_{(0,t]} \eta_p(\E[X\sm^p])\dd A\s + \alpha_2^p\EO[H^p\ti], \\
& \leq  (1-p) \int_{(0,t]} \eta_p(y(s^-))\dd A\s + \alpha_2^p\EO[H^p\ti],
\end{aligned}
\end{equation*}
which implies  by \cref{lemma:DetBihari}, $\tilde G_p(x) = (1-p)G(x^{1/p})$ and  $\tilde{G}_p^{-1}(x) = \big(G^{-1}\big(\tfrac{x}{1-p} \big)\big)^p$ (see \eqref{eq:G_umrechnen}):
\begin{equation*}
y(T) \leq  \tilde{G}_p^{-1}( \tilde{G}_p(\alpha_2^p \EO[H^p_T]) + (1-p)A_T) = G^{-1}(G(\alpha_2 \|H_T\|\pF) + A_T)^p
\end{equation*}
which implies by $\PO[X^*_T>u ] \leq \EO[X_{T\wedge \sigma_u}\wedge u] u^{-1}$
\begin{equation}
\begin{aligned}
\|X^*_T\|^p\pF & = \int_0^\infty  p \PP[X^*_T> u] u^{p-1} \dd u \\
& \leq p \int_0^\infty  \EO[X_{T\wedge \sigma_u}\wedge u] u^{p-2} \dd u  \\
& = \frac{1}{1-p} y(T)  \\
& \leq \frac{1}{1-p} G^{-1}(G(\alpha_2 \|H_T\|\pF ) + A_T))^p
\end{aligned}
\end{equation}
which implies the assertion \eqref{eq:claimConcaveNoSup-p}. \\

Now we prove \eqref{eq:claimConcaveNoSup-weak} for predictable $H$. Fix some $u>0, w>0$ and $t>0$ and choose $\lambda>0$ s.t. $\lambda u = w$. Due to predictability of $H$ there exists a localizing sequence $\tau^{(n)}_u$ for the stopping time $\tau_u \define \inf\{s \geq 0 \mid H\s \geq \lambda u\}$. Define $B\define {\{H_0 < \lambda u\}}$ and note that $\EO[H_{\tau^{(n)}_u \wedge t}\one_B] < \infty$. We will now use that we have shown \eqref{eq:claimConcaveNoSup-weak} already for the case that $\E[H_t]<\infty$ : This gives that  \eqref{eq:claimConcaveNoSup-weak} holds for $(X^{\tau^{(n)}_u}\one_B, A, H^{\tau^{(n)}_u}\one_B, M^{\tau^{(n)}_u}\one_B)$ and $\eta$, i.e.
\begin{equation*}
\PO[X^*_{t\wedge\tau^{(n)}_u}>u , H_0 < \lambda u] \leq \frac{1}{u} G^{-1}(G(\EO[H_{t\wedge\tau^{(n)}_u}]\one_B) + A_t) \leq \frac{1}{u} G^{-1}(G(\EO[H_{t} \wedge (\lambda u)]) + A_t).
\end{equation*}
Noting that $\{H_t < \lambda u\} \subseteq \{H_0 < \lambda u\}$ and that on $\{H_t < \lambda u\}$ we have $\lim_{n\to\infty} \tau_u^{(n)} = \tau_u > t$  gives:
\begin{equation*}
\begin{aligned}
\PO[X^*_t >u ] & \leq \sup_{n\in\N}\PO[X^*_{t\wedge\tau^{(n)}_u}>u , H_t < \lambda u] + \PO[H_t \geq \lambda u] 
\\
&\leq  
\frac{1}{u} G^{-1}(G(\EO[H_{t} \wedge (\lambda u)]) + A_t) + \PO[H_t \geq \lambda u].
\end{aligned}
\end{equation*}
\end{proof}

\subsection{Proof for concave \texorpdfstring{$\eta$}{eta} and random \texorpdfstring{$A$}{A}
(Corollary 3.8, Corollary 3.9)}

\begin{proof}[Proof of \cref{cor:randomA}]
We prove the claim for the case of \nameref{def:sup} and predictable $H$. The other cases follow by the same argument. We prove the asssertion by a time change \cite[Lemma 4.8]{GeissConvex}. Fix some $T>0$. We may assume w.l.o.g. that the continuous part of $A$ is strictly increasing and $A^c_\infty = \infty$, for details see \cite[Remark 4.10]{GeissConvex}. Let $(\tilde{X}, \tilde{A}, \tilde{H}, \tilde{M})$ denote the family of processes we obtain by applying \ \cite[Lemma 6.8]{GeissConvex} to 
$((X_{t\wedge T})_{t\geq 0}, A, (H_{t\wedge T})_{t\geq 0}, (M_{t\wedge T})_{t\geq 0})$. Due to $\tilde{A}_t = t$ being deterministic, we may apply \cref{thm:stochBihariConcave}, yielding
\begin{equation*}
\EO[(\tilde{X}_t^*)^p] \leq G^{-1}(G(\alpha_1 \alpha_2 \|\tilde{H}_t\|\pF) + \beta t)^p \qquad \forall t\geq 0.
\end{equation*}
Note that for any $t\geq 0$ we have
\begin{equation*}
\EO[\one_{\{A_T\leq t\}} (X^*_T)^p] \leq \EO[\sup_{r: A_r\leq t} (X_{r\wedge T})^p] = \EO[\sup_{r: A_r\leq t} (\tilde{X}_{A_r})^p] \leq \EO[\sup_{s\leq t} (\tilde{X}_s)^p] = \EO[(\tilde{X}_t^*)^p].
\end{equation*}
To simplify the notation in the following calculations, define 
\begin{equation*}
f:[0,\infty) \to [\alpha_1 \alpha_2 \|H_T\|\pF,\infty), \qquad x \mapsto G^{-1}(G(\alpha_1 \alpha_2 \|H_T\|\pF)+ \beta x) 
\end{equation*}
noting that concavity of $\eta$ implies that 
$[G(\alpha_1 \alpha_2 \|H_T\|\pF),\infty) \subseteq \textrm{domain}(G^{-1})$, hence implying that the map is well-defined. Moreover, $f$ is strictly increasing and bijective.
The two inequalities above yield together
\begin{equation}\label{eq:boundTimeShift}
\EO[\one_{\{A_T\leq t\}} (X^*_T)^p]\leq G^{-1}(G(\alpha_1 \alpha_2 \|H_T\|\pF) + \beta t)^p = f(t)^p
\end{equation}
  We have for $0<p<q$, $p<1$:
\begin{equation*}
\begin{aligned}
\EO[f(A_T)^{-q} (X_T^*)^p] & =\EO[q\int_{f(A_T)}^\infty x^{-q-1} dx (X_T^*)^p] \\
& =  q \int_{ \alpha_1 \alpha_2 \|H_T\|\pF }^\infty \EO[\one_{\{f(A_T) \leq  x\}} (X_T^*)^p]x^{-q-1}dx \\
& \overset{\eqref{eq:boundTimeShift}}{\leq}  q \int_{\alpha_1 \alpha_2 \|H_T\|\pF  }^\infty f(f^{-1}(x))^p  x^{-q-1}dx \\
& = q \int_{\alpha_1 \alpha_2 \|H_T\|\pF  }^\infty x^p x^{-q-1}dx \\
& \leq \frac{q}{q-p} (\alpha_1 \alpha_2 \|H_T\|\pF)^{p-q}.
\end{aligned}
\end{equation*}
This implies the first assertion. We obtain the upper bound for $\|X^*_t\|_{p/2,\F_0}$
by rewriting 
\begin{equation*}
\EO[(X_t^*)^{p/2}] = \EO\left [ (X_t^*)^{p/2} \big(G^{-1}(G(\alpha_1 \alpha_2 \|H_t\|\pF) + \beta A_t)\big)^{-q/2} \cdot \big(G^{-1}(G(\alpha_1 \alpha_2 \|H_t\|\pF) + \beta A_t)\big)^{q/2}\right].
\end{equation*}
and then applying H\"older's inequality.
\end{proof}

\begin{proof}[Proof of \cref{cor:randomAnoSup}]
We prove the claim for the case of $\E[H_T]<\infty$, the assertion can be shown by the same argument for predictable $H$. As in the proof of  \cref{cor:randomA} we prove the asssertion by a time change \cite[Lemma 4.8]{GeissConvex}. Fix some $T>0$. We may assume w.l.o.g. that the continuous part of $A$ is strictly increasing and $A^c_\infty = \infty$, for details see \cite[Remark 4.10]{GeissConvex}. Let $(\tilde{X}, \tilde{A}, \tilde{H}, \tilde{M})$ denote the family of processes we obtain by applying \ \cite[Lemma 6.8]{GeissConvex} to 
$((X_{t\wedge T})_{t\geq 0}, A, (H_{t\wedge T})_{t\geq 0}, (M_{t\wedge T})_{t\geq 0})$. Due to $\tilde{A}_t = t$ being deterministic, we may apply \cref{thm:stochBihariConcave}  $(\tilde{X}, \tilde{A}, \tilde{H}, \tilde{M})$, yielding for all $R>0, u>0$
\begin{equation*}
\PO[\tilde{X}^*_R > u] \leq \frac{G^{-1}(G(\EO[\tilde{H}_R]) + \tilde{A}_R)}{u}
\leq \frac{G^{-1}(G(\EO[H_T]) + R)}{u}
\end{equation*}
This implies using  $\{X^*_T > u, A_T \leq R\} \subseteq \{\tilde{X}^*_{A_T} > u, A_T \leq R\}  \subseteq  \{\tilde{X}^*_{R} > u\} $ 
\begin{equation*}
\PO[X^*_T > u ] \leq \PO[X^*_T > u, A_T \leq R] + \PO[A_T > R]  \leq \frac{G^{-1}(G(\EO[H_T]) + R)}{u} + \PO[A_T > R].
\end{equation*}
\end{proof}

%%%%%%%%%%%%%%%%%%%%%%%%%%%%%%%%%%%%%%%%%%%%%%%%%%%%%%%%%%%%%%%%%%
%%%%%%%%%%%%%%%%%%%%%%%%%%%%%%%%%%%%%%%%%%%%%%%%%%%%%%%%%%%%%%%%%%%
%%%%%%%%%%%%%%%%%%%%%%%%%%%%%%%%%%%%%%%%%%%%%%%%%%%%%%%%%%%%%%%%%%%

\subsection{Proof for general $\eta$ (\cref{thm:BihariRandomA})}
We prove  \cref{thm:BihariRandomA}\ref{item:SB-1}  by combining \cref{lemma:1} with \cref{rmk:Burkholder},
assertions \ref{item:SB-2}, \ref{item:SB-3} and \ref{item:SB-4} follow from \ref{item:SB-1}.

\begin{proof}[Proof of \cref{thm:BihariRandomA}]
\textbf{Proof of \ref{item:SB-1}:} Since for any stopping time $\tau$ the family of processes\\ $(X_{\cdot \wedge \tau}, A_{\cdot \wedge \tau}, H_{\cdot \wedge \tau},M_{\cdot \wedge \tau}) $ also satisfies assumption $\Asup$, it suffices to prove the claim for $\tau = T$. We may assume that $\E[H_T]<\infty$  and $\E[A_T] <\infty$ since otherwise the assertion is trivial. \cref{lemma:1} (using Assumption $\Asup$ and $E[H_T]<\infty$) implies for any $t\in[0,T]$
\begin{equation}\label{eq:proof:randomBihari1}
\E[X_t\wedge u] \leq \E[X_{t\wedge \sigma_u}\wedge u] \leq \E[\int_{(0,t]} \eta(X^*\sm)\one_{\{X^*\sm \leq u\}}\dd A\s] +  \E[H_T] \wedge u.
\end{equation}
Here, we used that by Fatou's lemma $\E[M^{L,u}_T] \leq 0$.  Set $F(x) \define G(x+\varepsilon) - G(\varepsilon)$ which is concave and satisfies $F(0)=0$ and $F'(\infty) = 0$. Combining \cref{rmk:Burkholder} with \eqref{eq:proof:randomBihari1} implies by Fubini's theorem:
\begin{equation*}
\begin{aligned}
\E[F(X_T)] &  = - \E\bigg[\int_{[0,\infty)} X_T \wedge u \dd F'(u) \bigg] \leq - \E\bigg[\int_{(0,T]}\eta(X^*\sm) \int_{X^*\sm}^\infty \dd F'(u) \dd A\s\bigg] + 
 F(\E[H_T]) \\
& \leq \E[A_T] + F(\E[H_T]).
\end{aligned}
\end{equation*}
Hence, this implies for all $\varepsilon>0$
\begin{equation}\label{eq:itemSB1}
\E[G(X_T+\varepsilon)] \leq \E[A_T] + G(\E[H_T]+\varepsilon).
\end{equation}
Letting $\varepsilon \to 0$ gives by monotone convergence
\begin{equation*}
\E[G(X_T)]  \leq \E[A_T] + G(\E[H_T]).
\end{equation*}

\textbf{Proof of \ref{item:SB-3}:}
Let $\delta>0$ be arbitrary. Set $\tau^{(n)}\define \inf\{t\geq 0 \mid X^{(n)}_t > \delta\}\wedge T$ using the definition $\inf \emptyset \define \infty$. Moreover, set 
$P^{(n)}_\delta  \define \PP[\sup_{t\in[0,T]}X_t^{(n)} >\delta]$. Recall that $X$ is right-continuous and $G$ is non-decreasing. Furthermore, in the proof of \ref{item:SB-1}
(see \eqref{eq:itemSB1}) 
we also proved $\E[G(X_{\tau\wedge T}+\varepsilon)] \leq \E[A_{\tau\wedge T}] + G(\E[H_{\tau\wedge T}]+\varepsilon)$ for all stopping times $\tau$ and all $\varepsilon>0$. Combining this yields:
\begin{equation*}
P^{(n)}_\delta G(\varepsilon + \delta) + (1-P^{(n)}_\delta)
G(\varepsilon)  \leq \E[G(X^{(n)}_{\tau^{(n)}\wedge T}+\varepsilon)] 
\leq \E[A_T] + G(E[H^{(n)}_T]+\varepsilon).
\end{equation*}
Rearranging the terms implies by continuity of $G$ for any $\varepsilon>0$ 
\begin{equation*}
\begin{aligned}
\limsup_{n\to\infty} P{(n)}_\delta \leq  \limsup_{n\to\infty}
\frac{\E[A_T] + G(E[H^{(n)}_T]+\varepsilon)- G(\varepsilon)}{G(\varepsilon +\delta) -G(\varepsilon)}
 = \frac{\E[A_T]}{G(\varepsilon +\delta) -G(\varepsilon)}.
\end{aligned}
\end{equation*}
Taking the limit $\varepsilon \to 0$ and using the assumption $\lim_{\varepsilon\to 0} G(\varepsilon) = -\infty$ implies the claim $\limsup_{n\to\infty} P^{(n)}_\delta=0$.

\textbf{Proof of \ref{item:SB-2}:} This follows immediately from \ref{item:SB-3} by choosing 
$(X^{(n)}, A, H^{(n)}, M^{(n)}) \define (X, A, H, M)$ for all $n$.

\textbf{Proof of \ref{item:SB-4}:} By assumption we have $X_t\geq c$ and $G(x) =\int_c^x \frac{\dd u}{\eta(u)}$ for all $x>0$, which implies that $G(X_t)\geq 0$ for all $t$. Hence,
applying Lenglart's inequality (see e.g. \cite[Corollaire II]{Lenglart} or \cite[Lemma 2.2 (ii)]{MehriScheutzow}) to assertion \ref{item:SB-1} yields $\|\sup_{t\in[0,T]} G(X_t)\|_p \leq \alpha_1 \alpha_2 \|A_T + G(\E[H_T])\|_p $.
\end{proof}

\subsection{Proof of the application (\cref{cor:solutions})}
\begin{proof}[Proof of \cref{cor:solutions}]
We may assume w.l.o.g. that $\lim_{x\to\infty}\eta_1(x) = +\infty$ and $\lim_{x\to\infty}\eta_2(x) = +\infty$. \\
\noindent \textbf{Proof of uniqueness of global solutions:} Let $(X_t)_{t\in[0,T]}$ and $(Y_t)_{t\in[0,T]}$  be two strong solutions of \eqref{eq:SDE} up to some deterministic time $T>0$ with the same initial condition $(z_t)_{t\in[-r,0]}$. For the uniqueness of global strong solutions it suffices to prove $\PP[\sup_{t\in[0,T]}|X_{t} - Y_{t}|^2 =0] = 1$. \\

\noindent Define
\begin{equation*}
\tau(R) \define \inf\{t\in[0,T] \mid |X_t|>R \text{ or } |Y_t| > R\}
\end{equation*}
where we set $\inf \emptyset \define T$. Then, by It\^o's formula we have for $t\in[0,T]$:
\begin{equation}\label{eq:proofUniqueness}
\begin{aligned}
|X_{t\wedge \tau(R)} - Y_{t\wedge \tau(R)}|^2 & = \int_0^{t\wedge \tau(R)} 2 \langle X_{s^-} - Y_{s^-}, f(s, X)-f(s, Y)\rangle \dd s \\
& + \int_0^{t\wedge \tau(R)} |h(s, X)-h(s, Y)|_F^2 \dd s \\
& +  \int_{(0, t\wedge \tau(R)]}\int_U |g(s, X,\xi) -  g(s, Y,\xi)|^2 \nu(\dd \xi) \dd s + M_{t\wedge \tau(R)} \\
 & \overset{(C1)}{\leq} \int_0^{t\wedge \tau(R)} L^R_s \eta_1 (\sup_{u\in[0,s]}|X_u-Y_u|^2 )\dd s + M_{t\wedge \tau(R)} \\
  & \leq \int_0^{t} L^R_s \eta_1(\sup_{u\in[0,s]}|X_{u\wedge \tau(R)}-Y_{u\wedge \tau(R)}|^2 )\dd s + M_{t\wedge \tau(R)},
\end{aligned}
\end{equation}
where the local martingale $M$ is given by
\begin{equation*}
\begin{aligned}
M\ti & \define \int_0^t 2\langle X_{s^-} - Y_{s^-}, (h(s, X)-h(s, Y)) \dd B\s \rangle  \\
& + \int_{(0,t]}\int_U 2\langle X_{s^-} - Y_{s^-}, g(s, X, \xi)-g(s, Y, \xi)\rangle \tilde{N}(\dd s, \dd \xi) \\
& + \int_{(0,t]} \int_{U} |g(s, X, \xi)-g(s, Y,\xi)|^2 \tilde{N}(\dd s, \dd \xi).
\end{aligned}
\end{equation*}
Applying \cref{thm:BihariRandomA}\ref{item:SB-2} to 
$(|X_{\cdot \wedge \tau(R)} - Y_{\cdot \wedge \tau(R)}|^2, 
\int_0^\cdot L_s^R \dd s, 0, M)$  and $\eta_1$ implies  $\PP[\sup_{t\in[0,T]}|X_{t\wedge \tau(R)} - Y_{t\wedge \tau(R)}|^2 =0] = 1$ for all $T> 0$, $R>0$ which yields the claim.

\bigskip

\noindent \textbf{Proof of the existence of strong global solutions:}
We show the existence of a solution using the Euler method, see e.g. \cite{Roeckner}, \cite{RenesseScheutzow}, \cite{MehriScheutzow}. We define for $n\in\N$ the Euler approximates as
\begin{equation}\label{eq:EulerApprox}
\begin{aligned}
X\nn_t & \define z_t \qquad \text{for } t\in[-r,0] \\
X_t\nn & \define X\nn_{\kn} + \int_{(\kn,t]}  f(s, X\nn_{\cdot \wedge \kn} ) \dd s
+ \int_{(\kn,t]} \int_U g(s, X\nn_{\cdot \wedge\kn}, \xi) \tilde{N}(\dd s, \dd \xi)  \\
& \qquad + \int_{(\kn,t]} h(s,  X\nn_{\cdot \wedge \kn}) \dd B\s \qquad
\text{for } t\in\big(\tfrac{k}{n}, \tfrac{k+1}{n}\big], k\in\N_0
\end{aligned}
\end{equation}
The Euler approximate $X\nn$ has c\`adl\`ag paths and is adapted. In particular, the stochastic integrals are well-defined. Define $k(n,t) \define \kn$ for $t\in(\kn, \frac{k+1}{n}]$, $k\in\N_0$ and $k(n,t)\define t$ for $t\in[-r,0]$. 
\bigskip
\noindent With the notation $k(n,t)$ \eqref{eq:EulerApprox} can be rewritten as
\begin{equation}\label{eq:EulerApprox2}
\begin{aligned}
X_t\nn &= z_0 + \int_{0}^t f(s, X\nn_{\cdot \wedge k(n,s)} ) \dd s
+ \int_{(0,t]} \int_U g(s, X\nn_{\cdot \wedge k(n,s)}, \xi) \tilde{N}(\dd s, \dd \xi)  \\
& \qquad + \int_0^t h(s, X\nn_{\cdot \wedge k(n,s)}) \dd B\s \qquad \text{for all } t\geq 0.
\end{aligned}
\end{equation}
Define the remainder 
\begin{equation*}
p_t\nn \define X\nn_{k(n,t)} - X_t\nn, \qquad t\in[-r,\infty).
\end{equation*}
The stochastic process $p\nn$ is adapted and $p\nn_{k/n^+} = 0$ for every $k\in\N_0$.  Note that $t\mapsto X\nn_{k(n,t)},\, t\geq 0$ is c\`agl\`ad and $X\nn_t$ is \cadlagg. Therefore, $p_t\nn$ is neither c\`agl\`ad nor c\`adl\`ag.
Fix $T\geq 0$ and define for $R>0$ the stopping times 
\begin{equation*}
\tau\nn_R \define \big(\inf\big\{t\geq 0 \mid |X\nn_t| > \frac{R}{3} \big\} \wedge T\big)\one_{\{R>3\sup_{s\in[-r,0]} |z\s|\}}.
\end{equation*}
Then
\begin{equation}\label{eq:pnBound}
|p\nn_t| \leq \frac{2R}{3}, \quad |X\nn_t| \leq \frac{R}{3}, \quad t\in(0,\tau_R\nn).
\end{equation}
On $\{R>3\sup_{s\in[-r,0]} |z\s|\}$ the inequalities extend to $[-r,\tau_R\nn)$ and we have $\tau_R\nn>0$ (due to the right-continuity of $X\nn$).

\bigskip

As $T$ was arbitrary, it suffices to prove the existence of a strong solution on  $[0,T]$. The existence proof is done in the following steps:
\begin{enumerate}
\item\label{item:stepi} For every $t \geq 0$, we have $\one_{(0,\tau_R^{(n)})}(t) \sup_{u\in (k(n,t),t]}|p_u^{(n)}| \to 0$ in probability as $n\to\infty$.
\item\label{item:stepii} 
$\bigg\|G_{2}(1+\sup_{u\in[0,T]} |X\nn_{u\wedge \tau_R^{(n)}}|^2)\bigg\|_{1/2} \leq C(T,R,n)$  for some $C(T,R,n)$ satisfying
\begin{equation*}
\lim_{n\to\infty} C(T,R,n) = \tilde{C}(T) \text{ for all } R>0.
\end{equation*}
\item\label{item:stepiii} $\lim_{R\to\infty} \limsup_{n\to\infty} \PP[\tau_R^{(n)} < T]= 0$.
\item\label{item:stepiv} $\lim_{n,m\to\infty} \PP\big[ \sup_{t\in[0,T]} |X_t^{(n)}-X_t^{(m)}| > \varepsilon]=0$ for all $\varepsilon>0$.
\item\label{item:stepv} There exists a \cadlagg adapted process $(X_t)_{t\geq 0}$ such that \\$\lim_{n\to\infty} \PP[\sup_{t\in[0,T]} |X_t^{(n)} - X_t| > \varepsilon] = 0$  for all $\varepsilon>0$ and $X$ is a strong solution of equation \eqref{eq:SDE} on $[0,T]$.
\end{enumerate}
We only prove the steps \ref{item:stepii}, \ref{item:stepiii}, \ref{item:stepiv}: The remaining steps are proven as in \cite[Theorem 3.3]{MehriScheutzow}.

%Proof of \ref{item:stepi}: This step is proven using \cref{item:c4}, the proof is identical to \cite{MehriScheutzow}. \\

\noindent \textbf{Proof of \ref{item:stepii}:}  Using It\^o's formula, we have for $t\in[0,T]$:
\begin{equation*}
\begin{aligned}
|X^{(n)}_t|^2 & = |z_0|^2 + \int_0^t 2\langle X^{(n)}_{s^-}, f(s,  X\nn_{\cdot \wedge k(n,s)})\rangle \dd s  + \int_{(0,t]}\int_{U} |g(s,  X\nn_{\cdot \wedge k(n,s)}, \xi)|^2 \nu (\dd \xi) \dd s \\
&  \quad  + \int_0^t |h(s, X\nn_{\cdot \wedge k(n,s)})|_F^2 \dd s + M_t^{(n)},
\end{aligned}
\end{equation*}
where $M\nn$ is a local martingale defined by
\begin{equation*}
\begin{aligned}
M_t\nn & \define \int_{(0,t]} \int_U 2 \langle \xs, \g \rangle \tilde{N}(\dd s,\dd \xi) \\
& \quad + \int_{(0,t]} \int_U |\g|^2 \tilde{N}(\dd s, \dd \xi) \\
& \quad + \int_0^t 2\langle \xs, \h \dd B\s\rangle.
\end{aligned}
\end{equation*}
Using (C2) and (C4), we have (using $\lim_{r\nearrow s} X\nn_{r\wedge k(n,s)} =  X\nn_{k(n,s)} = \lim_{r\nearrow s} X\nn_{k(n,r)}$ for $s>0$):
\begin{equation}\label{eq:stepii-13}
\begin{aligned}
|X^{(n)}_{t\wedge\tau_R^{(n)}}|^2 & \overset{(C2)}{\leq} |z_0|^2 + \int_0^{t\wedge\tau_R^{(n)}} 2 \langle \xs- \lim_{r\nearrow s} X\nn_{r\wedge k(n,s)}, f(s, X\nn_{\cdot \wedge k(n,s)})\rangle \dd s \\
& \quad + \int_0^{t\wedge\tau_R^{(n)}} K_s \, \eta_2 \big(1 + \sup_{u\in[-r,s]} |X^{(n)}_{u\wedge k(n,s)}|\big)\dd s + M_t^{(n)} \\
&  \overset{(C4)}{\leq}  |z_0|^2 + \int_0^{t\wedge\tau_R^{(n)}} 2 |p^{(n)}_{s^-}|\tilde{K}^R_s \dd s \\
& \quad + \int_0^t K_s \, \eta_2 \big(1 + \sup_{u\in[-r,0]} |z_u|^2 + \sup_{u\in[0,s]} |X^{(n)}_{u\wedge\tau_R^{(n)}} |^2 \big) \dd s + M_{t\wedge\tau_R^{(n)}}^{(n)}
\end{aligned}
\end{equation}
Set
\begin{equation*} 
H_t^{n,R}\define 1 + \sup_{u\in[-r,0]} |z_u|^2 + |z_0|^2 + \int_0^{t\wedge\tau_R^{(n)}} 2 |p^{(n)}_{s^-}|\tilde{K}_R(s) \dd s.
\end{equation*}
For almost every $s\in[0,T]$ we have \eqref{eq:pnBound}, \ref{item:stepi}  and dominated convergence
\begin{equation*}
0 \leq \limsup_{n\to\infty} \E[ |p_{s^-}^{(n)}| \one_{\{s < t\wedge\tau_R^{(n)}\}}\tilde{K}^R_s] 
\leq \limsup_{n\to\infty} \E[ \sup_{u\in(k(n,s),s]}|p_{u}^{(n)}| \one_{\{s < t\wedge\tau_R^{(n)}\}}\tilde{K}^R_s]  \overset{\ref{item:stepi}}{=} 0.
\end{equation*}
Therefore, by dominated convergence, we have for all $t\in[0,T]$:
\begin{equation*}
\lim_{n\to\infty} \E\bigg[ \int_0^{t\wedge\tau_R^{(n)}} 2 |p^{(n)}_{s^-}|\tilde{K}^R_s \dd s\bigg]  = \lim_{n\to\infty}  \int_0^{t} 2 \E[ |p^{(n)}_{s^-}| \one_{\{s < t\wedge\tau_R^{(n)}\}}\tilde{K}^R_s]\dd s  = 0.
\end{equation*}
In particular, by (C5) this implies
\begin{equation}\label{eq:stepii-1}
\lim_{n\to\infty} \E[H^{n,R}_T]  = 1 + \E[\sup_{u\in[-r,0]} |z_u|^2] + \E[|z_0|^2] <\infty.
\end{equation}
We apply \cref{thm:BihariRandomA}\ref{item:SB-4} to $Y\ti \define 1 + \sup_{u\in[-r,0]} |z_u|^2 + |X^{(n)}_{t\wedge\tau_R^{(n)}\wedge T}|^2$, $t\geq  0$, which satisfies due to \eqref{eq:stepii-13} $Y\ti \leq \int_0^t K_s \eta_2(Y^*\sm) ds + M_{t\wedge\tau_R^{(n)}\wedge T} \nn + H^{n,R}_{t\wedge\tau_R^{(n)}\wedge T}$. 
This yields (noting that $\alpha_1 \alpha_2 = 8$ for $p=1/2$):
\begin{equation}\label{eq:stepii-2}
\begin{aligned}
\bigg\|G_2(1+\sup_{t\in[0,T]}|X\nn_{t\wedge\tau_R^{(n)}}|^2)\bigg\|_{1/2}
& \leq \bigg\|G_2(Y^*_T)\bigg\|_{1/2} \\
& \leq 8\bigg\|\int_0^T K_t \dd t + G_2(\E[H^{n,R}_T])\bigg\|_{1/2} =\vcentcolon C(T,R,n).
\end{aligned}
\end{equation}
Combining \eqref{eq:stepii-1} with \eqref{eq:stepii-2} implies by continuity of $G_2$ for 
\begin{equation*}
\tilde{C}(T) \define 8 \bigg\|\int_0^T K_t \dd t + G_2(1 + \E[\sup_{u\in[-r,0]} |z_u|^2] + \E[|z_0|^2])\bigg\|_{1/2} 
\end{equation*}
the assertion. \\

\noindent \textbf{Proof of \ref{item:stepiii}:} We have due to $G_2$ being non-decreasing, $G_2(x) \geq 0$ for $x\geq 1$  and $\lim_{x\to\infty} G_2(x) =+\infty$:
\begin{equation*}
\begin{aligned}
\limsup_{R\to\infty} \limsup_{n\to\infty}
\PP[\tau_R\nn < T] & 
\leq \limsup_{R\to\infty} \limsup_{n\to\infty} \PP[\sup_{t\in[0,\tau\nn_R\wedge T]} |X\nn_t| \geq R/3]
+ \lim_{R\to\infty} \PP[\sup_{s\in[-r,0]} |z\s|\geq R/3]\\
& \overset{(C5)}{\leq}  \limsup_{R\to\infty} \limsup_{n\to\infty}\frac{ \E[G_2(1+\sup_{t\in[0,\tau^{(n)}_R\wedge T]} |X\nn_t|^2)^{1/2}]}{G_2(1+(R/3)^2)^{1/2}} \\
&  \overset{\ref{item:stepii}}{\leq} \limsup_{R\to\infty} \limsup_{n\to\infty} \frac{C(T,R,n)^{1/2}}{G_2(1+(R/3)^2)^{1/2}} \\
& \overset{\ref{item:stepii}}{\leq} \limsup_{R\to\infty} \frac{\tilde{C}(T)^{1/2}}{G_2(1+(R/3)^2)^{1/2}} \\
& = 0.
\end{aligned}
\end{equation*}

\noindent \textbf{Proof of \ref{item:stepiv}:} 
Set $\tRnm \define \tau\nn_R \wedge \tau\mm_R$. It\^o's formula gives for $t\in[0,T]$:
\begin{equation*}
\begin{aligned}
| X_t^{(n)} - X_t^{(m)}|^2 & = \int_0^t 2 \langle X\nn_{s^-} - X\mm_{s^-}, f(s, X\nn_{\cdot \wedge k(n,s)})- f(s, X\mm_{\cdot \wedge k(m,s)}) \rangle \dd s 
\\ & 
+ \int_{(0,t]} \int_{U} |g(s, X\nn_{\cdot \wedge k(n,s)}, \xi)- g(s,  X\mm_{\cdot \wedge k(m,s)}, \xi) |^2 \nu (\dd \xi) \dd s
\\ & 
+ \int_0^t |h(s, X\nn_{\cdot \wedge k(n,s)})- h(s, X\mm_{\cdot \wedge k(m,s)})|^2_F \dd s + M_t^{n,m} 
\end{aligned}
\end{equation*}
where $(M_t^{n,m})_{t\geq 0}$ is a local martingale starting in $0$. Using (C1) we obtain:
\begin{equation}\label{eq:proofExistence2}
\begin{aligned}
|X_{t\wedge \tRnm}\nn - X_{t\wedge\tRnm}\mm|^2 \overset{(C1)}{\leq} 
\int_0^{t\wedge\tRnm} L^R_s \eta_{1}\bigg( \sup_{u\in[-r,s]} |X\nn_{u\wedge k(n,s)} -X\mm_{u\wedge k(m,s)}|^2\bigg)\dd s \\
- 2\int_0^{t\wedge\tRnm}\langle p\nn_{s^-} - p\mm_{s^-}, f(s, X\nn_{\cdot \wedge k(n,s)})- f(s, X\mm_{\cdot \wedge k(m,s)}) \rangle \dd s  +  M_{t\wedge\tRnm}^{n,m}.
\end{aligned}
\end{equation}
To be able to apply  \cref{thm:BihariRandomA}, we need to  compute $\sup_{u\in[-r,s]} |X\nn_{u\wedge k(n,s)} -X\mm_{u\wedge k(m,s)}|^2$. We have for $u\in(0,s]$:
\begin{equation*}
\begin{aligned}
X_{u\wedge k(n,s)} &= \begin{cases}
 X\nn_{u} & \text{ for } u\in[0,k(n,s)] \\
X\nn_{k(n,s)} =  X\nn_{k(n,u)} = p\nn_u + X\nn_{u}
& \text{ for } u\in(k(n,s),s] 
\end{cases} \\[1em]
& =  X\nn_{u} +  p\nn_u \one_{\{ u \in(k(n,s),s] \}},
\end{aligned}
\end{equation*}
which yields
\begin{equation*}
|X\nn_{u\wedge k(n,s)} -X\mm_{u\wedge k(m,s)}|^2
 \leq 2|X\nn_{u} -X\mm_{u}|^2 + 4 |p\nn_u|^2 \one_{\{ u \in(k(n,s),s] \}} + 4|p\mm_u|^2 \one_{\{ u \in(k(m,s),s] \}}.
\end{equation*}
Using the notation
\begin{equation*}
\begin{aligned}
p^{n,m,R}_s &\define  4\one_{(0,\tau\nn_R)}(s)\sup_{u\in(k(n,s),s]}|p_u\nn|^2 + 4
\one_{(0,\tau\mm_R)}(s)\sup_{u\in(k(m,s),s]}|p_u\mm|^2, \\
H^{n,m,R}_t &\define \int_0^t \one_{(0,\tRnm)}(s) 4 (|p\nn_{s^-}|+|p\mm_{s^-}|)\tilde{K}^R_s \dd s \\ &\quad 
+ \int_0^{t\wedge\tRnm} L^R_s \eta_{1}\bigg( 2\sup_{u\in[-r,s]} |X\nn_{u\wedge\tRnm} -X\mm_{u\wedge\tRnm} |^2 + p^{n,m,R}_s \bigg) \\ & \quad - L^R_s \eta_{1}\bigg( 2\sup_{u\in[-r,s]} |X\nn_{u\wedge\tRnm} -X\mm_{u\wedge\tRnm} |^2  \bigg) \dd s,
\end{aligned}
\end{equation*}
we obtain from \eqref{eq:proofExistence2} using (C4)
\begin{equation}\label{eq:proofExistence3}
\begin{aligned}
|X_{t\wedge \tRnm}\nn - X_{t\wedge\tRnm}\mm|^2  & \overset{(C4)}{\leq} 
\int_0^{t} L^R_s \eta_{1}\bigg( 2\sup_{u\in[-r,s]} |X\nn_{u\wedge\tRnm} -X\mm_{u\wedge\tRnm} |^2 \bigg) \dd s \\
& \qquad + M_{t\wedge\tRnm}^{n,m}  + H^{n,m,R}_t.
\end{aligned}
\end{equation}
Note that due to \eqref{eq:pnBound} there exists a constant $\hat{C}(R)>0$ that only depends on $R$ such that  
\begin{equation*}
2\sup_{u\in[-r,s]} |X\nn_{u\wedge\tRnm} -X\mm_{u\wedge\tRnm} |^2 \one_{\{s< \tRnm\}} \leq \hat{C}(R),\qquad 
p^{n,m,R}_s \leq \hat{C}(R).
\end{equation*}
Hence, by uniform continuity of $\eta_1$ on $[0,2\hat{C}(R)]$ and \ref{item:stepi},  for any fixed 
$s\in[0,t]$, we have
\begin{equation*}
\eta_{1}\bigg( 2\sup_{u\in[-r,s]} |X\nn_{u\wedge\tRnm} -X\mm_{u\wedge\tRnm} |^2 + p^{n,m,R}_s \bigg) - \eta_{1}\bigg( 2\sup_{u\in[-r,s]} |X\nn_{u\wedge\tRnm} -X\mm_{u\wedge\tRnm} |^2  \bigg) \to 0  
\end{equation*}
in probability for $m,n\to\infty$. Hence, the integrand of $H^{n,m,R}_t$ converges to $0$ in probability for almost every $s\in[0,t]$. 
Therefore, by dominated convergence we have:
\begin{equation*}
\limsup_{n, m\to\infty} \E[H^{n,m,R}_T] = 0.
\end{equation*}
We multiply \eqref{eq:proofExistence3} by $2$  and apply 
\cref{thm:BihariRandomA}\ref{item:SB-3}  to
\begin{equation*}
(2 \sup_{s\in[0,t]}|X_{s\wedge \tRnm}\nn - X_{s\wedge\tRnm}\mm|^2, \,\,  2 \int_0^t L^R_s \dd s, \,\,  2 H^{n,m,R}_t,  \,\,  2 M_{t\wedge\tRnm}^{n,m}).
\end{equation*} 
(This can be done by using that $\lim_{n,m\to\infty} a_{n,m} = a$ if and only if for all $(n_k)_k, \,(m_k)_k \subseteq \N$ with 
$\lim_{k\to\infty} n_k = \lim_{k\to\infty} m_k = \infty$
 we have $\lim_{k\to\infty} a_{n_k,m_k} = a$.)
We obtain that $\lim_{n,m\to\infty} \PP[\sup_{t\in[0,\tRnm\wedge T]} |X_t\nn -X_t\mm| > a]  =0$ for any $a>0$. 

By \ref{item:stepiii} we have:
\begin{equation*}
\begin{aligned}
& \limsup_{n,m\to\infty}
\PP[\sup_{t\in[0,T]} |X_t\nn - X\mm_t| > a] \\
& \quad \leq 
\limsup_{R\to\infty}\limsup_{n,m\to\infty}\big(\PP[T > \tau_R\nn] + \PP[T>\tau_R\mm] + \PP[\sup_{t\in[0,\tRnm\wedge T]} |X_t\nn -X_t\mm|>a])\\
& \quad =0,
\end{aligned}
\end{equation*}
which implies the claim.
\end{proof}

\section{Appendix}
\subsection{Proof of the (deterministic) Bihari-LaSalle inequality}

\begin{proof}[Proof of \cref{lemma:DetBihari}]
Denote by $y(t)$ the right-hand side of \eqref{eq:detBihari}, which is a non-decreasing right-continuous function. Define $t_{k,n}\define  k t 2^{-n}$ for $k,n\in\N$. We obtain, using that $\eta$ is non-decreasing
\begin{equation*}
\begin{aligned}
G(y(t))- G(y(0)) & =  \lim_{n\to\infty} \sum_{k=1}^{2^n} \big[ G(y(t_{k,n})) - G(y(t_{k-1,n})) \big] \\
& =  \lim_{n\to\infty} \sum_{k=1}^{2^n} \int_{y( t_{k-1,n})}^{y(t_{k,n})} \frac{\dd u}{\eta(u)} \\
& \leq   \lim_{n\to\infty} \sum_{k=1}^{2^n} \frac{1}{\eta(y(t_{k-1,n}))}
\big[ (y(t_{k,n}) - y(t_{k-1,n})\big] \\
& \overset{*}{=} \int_{(0,t]} \frac{\dd y(s) }{\eta(y(s^-))}  
\leq \int_{(0,t]} \frac{\dd y(s) }{\eta(x(s^-))}  
= A(t).
\end{aligned}
\end{equation*}
We use the left-continuity of $\eta$ in the equality $(*)$. Rearranging the terms, noting that $y(0)=H$ and applying $G^{-1}$ (if possible) implies the assertion.
\end{proof}

\subsection{Counterexample: Random integrators $A$ in Theorem 3.1}
Let $(X,A,H,M)$  and $\eta$ satisfy \nameref{def:nosup}. The following counterexample shows: For random integrators $A$ we may not expect in general estimates of the type
\begin{equation*}
\|X^*_T\|_p  \leq c_1 \big\|G^{-1}\big(G\big(c_2 H_T \big) + A_T\big)\big\|_p
\end{equation*}
to hold, where $c_1$ and $c_2$ are constants that only depend on $p\in(0,1)$. This type of estimate fails to hold true even for constant $H$ and $\eta(x) \equiv x$.
\begin{example} \label{ex:counterexample2}
For some fixed $p\in(0,1)$ we will construct stochastic processes $X_\gamma = (X_{t,\gamma})_{t\geq 0}$, $A_\gamma=(A_{t,\gamma})_{t\geq 0}$ and $M_\gamma = (M_{t,\gamma})_{t\geq 0}$ depending on a parameter $\gamma\in[0,\infty)$, such that $(X_\gamma, A_\gamma, 1, M_\gamma)$ satisfy $\Anosup$, the equality 
\begin{equation}\label{eq:counterexample}
X_{t,\gamma} = \int_0^t X_{s^-,\gamma} \dd A_{s,\gamma} + M_{t,\gamma} + 1, \qquad  t\geq 0,
\end{equation}
and
\begin{equation*}
\lim_{n\to\infty} \frac{\|X^*_{T_n,\gamma_n}\|_p}{\|\e^{A_{T_n,\gamma_n}}\|_p}  = \infty
\end{equation*}
for suitably chosen sequences $(\gamma_n)_{n\in\N},(T_n)_{n\in\N}  \subset(0,\infty)$. \\

Choose $(\Omega, \F, \PP) \define ([0,1], \mathcal{L}[0,1], \dd s)$, $\F_t \define \mathcal{N}$ for $t\in[0,1)$ and $\F_t \define \F$ for $t\geq 1$, where $\mathcal{L}[0,1]$ denotes the Lebesgue $\sigma$-field on $[0,1]$ and $\mathcal{N}$ denotes the $\sigma$-field generated by the Lebesgue null sets. We define the martingale by
\begin{equation*}
\begin{aligned}
M_{t,\gamma}(\omega) & \define \one_{[1,\infty)}(t)(\gamma \one_{[0,\frac{\e}{\e + \gamma}]}(\omega) - \e\one_{(\frac{\e}{\e+\gamma},1]}(\omega)) \qquad \forall t\geq 0, \omega \in\Omega, \\
\tau_{\gamma}(\omega) & \define \infty \one_{[0,\frac{\e}{\e+\gamma}]}(\omega) + \one_{(\frac{\e}{\e+\gamma},1]}(\omega), \\
A_{t,\gamma} & \define t\wedge \tau_{\gamma}.
\end{aligned}
\end{equation*}
It can be easily verified, that 
\begin{equation*}
X_{t,\gamma}(\omega) \define 
\begin{cases} \e^t \qquad &\forall t\in[0,1),\, \omega\in\Omega  \\
 \e^{t-1}(e^1+\gamma) &\forall t\in[1,\infty),\, \omega \in[0,\frac{\e}{\e+\gamma}] \\
 0 &\forall t\in[1,\infty),\, \omega \in(\frac{\e}{\e+\gamma},1].
\end{cases}
\end{equation*}
satisfies \eqref{eq:counterexample}. This implies for $T\geq 1$:
\begin{equation*}
\begin{aligned}
\frac{\|X^*_{T,\gamma}\|_p^p}{
\|\exp(A_{T,\gamma})\|_p^p} & = \frac{ \e^{p(T-1)}(\e+\gamma)^p \PP\{[0,\frac{\e}{\e + \gamma}]\} + \e^p \PP\{ (\frac{\e}{\e+\gamma},1]\}}{ \e^{p T} \PP\{[0,\frac{\e}{\e + \gamma}]\} + \e^{p} \PP\{ (\frac{\e}{\e+\gamma},1]\}}  \\
 & = \frac{ \e^{-p}(\e+\gamma)^p + \e^p \PP\{ (\frac{\e}{\e+\gamma},1]\} \PP\{[0,\frac{\e}{\e + \gamma}]\}^{-1}\e^{-pT} }{1 + \e^{p} \PP\{ (\frac{\e}{\e+\gamma},1]\}\PP\{[0,\frac{\e}{\e + \gamma}]\}^{-1}\e^{-pT}}.  
\end{aligned}
\end{equation*}
Set $\gamma_n \define n$ for all $n\in\N$. For each $\gamma_n$ we can find a $T_n$ such that 
\begin{equation*} 
\frac{\|X^*_{T_n,\gamma_n}\|_p^p}{
\|\e^{A_{T_n,\gamma_n}}\|_p^p} \geq  \e^{-p}(\e+n)^p-1,
\end{equation*}
which implies the claim.
\end{example}

\renewcommand{\abstractname}{Acknowledgements}
\begin{abstract}
The author wishes to thank Michael Scheutzow for his valuable suggestions and comments, and Mark Veraar for his idea to include the weak $L^1$ estimates, i.e. \eqref{eq:claimConcaveNoSup-weak} in \cref{thm:stochBihariConcave} and
\cref{cor:randomAnoSup}.
\end{abstract}

%\bibliographystyle{plain} % We choose the "plain" reference style
%\bibliography{GronwallCalc2} 

\end{document}